\documentclass[11pt,a4paper]{article}

\usepackage{amsthm,amssymb,amsmath,amscd,latexsym}
\usepackage[all]{xy}

\theoremstyle{plain}
\newtheorem{thm}{Theorem}[section]
\newtheorem{lem}[thm]{Lemma}
\newtheorem{prop}[thm]{Proposition}
\newtheorem{cor}[thm]{Corollary}

\theoremstyle{definition}
\newtheorem{dfn}[thm]{Definition}
\newtheorem{ex}[thm]{Example}
\newtheorem{rem}[thm]{Remark}

\numberwithin{equation}{section}

\newcommand{\bA}{\mathbb{A}}

\newcommand{\bC}{\mathbb{C}}
\newcommand{\bF}{\mathbb{F}}

\newcommand{\bR}{\mathbb{R}}

\newcommand{\bZ}{\mathbb{Z}}

\newcommand{\cF}{\mathcal{F}}
\newcommand{\hcF}{\what{\cF}}

\newcommand{\cH}{\mathcal{H}}
\newcommand{\cK}{\mathcal{K}}
\newcommand{\cL}{\mathcal{L}}

\newcommand{\cO}{\mathcal{O}}
\newcommand{\cP}{\mathcal{P}}
\newcommand{\cR}{\mathcal{R}}
\newcommand{\cS}{\mathcal{S}}
\newcommand{\cU}{\mathcal{U}}

\newcommand{\bfF}{\mathbf{F}}

\newcommand{\bfR}{\mathbf{R}}
\newcommand{\bfV}{\mathbf{V}}

\newcommand{\rpb}{\mathcal{P}}
\newcommand{\rpbv}{\mathcal{P}^{\vee}}

\newcommand{\Sky}{\mathbf{Sky}}

\newcommand{\sky}{\Sky (\what T)}
\newcommand{\skynb}{\Sky_n (\what{T}_b)}
\newcommand{\skyb}{\Sky (\what{T}_b)}
\newcommand{\Loc}{\mathbf{Loc}}

\newcommand{\loct}{\Loc (T)}

\newcommand{\loctb}{\Loc (T_b)}
\newcommand{\vectn}{\mathbf{Vect}_n^\nabla (M)}

\newcommand{\vectno}{\mathbf{Vect}_n^\npart (M)}
\newcommand{\vecto}{\mathbf{Vect}^\npart (M)}
\newcommand{\specn}{\mathbf{Spec}_n^\nabla (\what M)}

\newcommand{\relskyn}{\mathbf{RelSky}_n (\what M)}
\newcommand{\relsky}{\mathbf{RelSky} (\what M)}

\newcommand{\locmn}{\mathbf{Loc}_n (M)}

\newcommand{\lochmn}{\mathbf{SpecLoc}_n (\what M)}
\newcommand{\Lochmn}[1]{\mathbf{SpecLoc}_{#1} (\what M)}

\newcommand{\Vect}[2]{\mathbf{Vect}_{#1}^\nabla (#2)}
\newcommand{\Vecto}[2]{\mathbf{Vect}_{#1}^\npart (#2)}
\newcommand{\Vct}[2]{\mathbf{Vect}_{#1} (#2)}
\newcommand{\Spec}[1]{\mathbf{Spec}_{#1}^\nabla (\what M)}

\newcommand{\Cinfty}[2]{C^{\infty}(#1, #2)}
\newcommand{\Morph}{\operatorname{Morph}}
\newcommand{\Id}{\operatorname{id}}

\newcommand{\id}{\mathbf{1}}

\newcommand{\what}{\widehat}
\newcommand{\wbar}{\overline}

\newcommand{\wti}{\widetilde}
\newcommand{\bsl}{\backslash}

\newcommand{\Adj}{\operatorname{Adj}}
\newcommand{\Aut}{\operatorname{Aut}}

\newcommand{\End}{\operatorname{End}}

\newcommand{\Hom}{\operatorname{Hom}}

\newcommand{\GL}{\operatorname{GL}}

\newcommand{\U}{\operatorname{U}}

\newcommand{\Ob}{\operatorname{Ob}}

\newcommand{\rank}{\operatorname{rank}}

\newcommand{\reg}{\operatorname{reg}}
\newcommand{\br}{\operatorname{sing}}
\newcommand{\supp}{\operatorname{supp}}

\newcommand{\Tr}{\operatorname{Tr}}

\newcommand{\Diff}{\operatorname{Diff}}

\newcommand{\fin}{\operatorname{fin}}

\newcommand{\ASD}{\operatorname{(A)SD}}

\newcommand{\ra}{\rightarrow}
\newcommand{\Ra}{\Rightarrow}
\newcommand{\la}{\leftarrow}
\newcommand{\lra}{\longrightarrow}

\newcommand{\hra}{\hookrightarrow}
\newcommand{\bra}{\mapsto}

\newcommand{\oset}[1]{\overset {#1}{\ra}}
\newcommand{\osetl}[1]{\overset {#1}{\lra}}
\newcommand{\ocong}[1]{\overset {#1}{\cong}}
\newcommand{\ot}{\otimes}
\newcommand{\abs}[1]{\vert {#1} \vert}

\newcommand{\fm}{\mathbf F}
\newcommand{\fmh}{\what{\fm}}
\newcommand{\ps}{p_{\Sigma}}
\newcommand{\psstar}{p_{\Sigma, *}}
\newcommand{\hps}{\hat{p}_{\Sigma}}
\newcommand{\hpsstar}{\hat{p}_{\Sigma,*}}
\newcommand{\hpis}{\hat{\pi}_{\Sigma}}
\newcommand{\trivcon}{\underline{\operatorname{d}}}
\newcommand{\Pv}{{P}^{\vee}}
\newcommand{\npart}{{\overset{{\operatorname{o}}}\nabla}}

\newcommand{\vphi}{\varphi}

\title{A Fourier--Mukai transform for real torus bundles}
\author{James F. Glazebrook \\ Eastern Illinois University \\
Department of Mathematics and Computer Science \\
Charleston, IL 61920 USA \\ cfjfg@eiu.edu \\ \\ 
Marcos Jardim \\University of Massachusetts at Amherst \\
Department of Mathematics and Statistics \\
Amherst, MA 01003--9305 USA \\ jardim@math.umass.edu \\ \\ Franz
W. Kamber \\ University of Illinois at Urbana--Champaign \\
Department of Mathematics \\ Urbana, IL 61801 USA\\
kamber@math.uiuc.edu}

\begin{document}
\maketitle

\begin{abstract}
We construct a Fourier--Mukai transform for smooth complex vector 
bundles $E$ over a torus bundle $\pi:M \to B~,$ the vector bundles 
being endowed with various structures of increasing complexity. 
At a minimum, we consider vector bundles $E$ with a flat partial 
unitary connection, that is families or deformations of flat 
vector bundles (or unitary local systems) on the torus $T~.$ 
This leads to a correspondence between such objects on $M$ and 
relative skyscraper sheaves $\cS$ supported on a spectral covering 
$\Sigma \hra \what M~,$ where $\hat\pi:\what{M} \to B$ is the flat 
dual fiber bundle.  Additional structures on $(E,\nabla)$ (flatness, anti-self-duality)
will be reflected by corresponding data on the transform $(\cS, \Sigma)~.$ 
Several variations of this construction will be presented, 
emphasizing the aspects of foliation theory which enter into 
this picture.


\end{abstract}

\vfill
\newpage

\tableofcontents

\vfill
\newpage
\baselineskip20pt


\section{Introduction} \label{intro}

The construction nowadays known as the {\em Fourier--Mukai transform} first
appeared in a seminal work by Mukai \cite{Mu}, where it was shown that
the derived categories of sheaves on an abelian variety (e.g. a complex
torus) is equivalent to the derived category of coherent sheaves on the
{\em dual} abelian variety. 

Since then, the Fourier--Mukai transform has been generalized in different ways, 
and has led to a number of interesting results concerning not only the derived
categories of coherent sheaves, but also the moduli spaces of stable sheaves on
abelian varieties, K3 surfaces and elliptic surfaces. 

This paper draws on two types of generalization of Mukai's original ideas. 
First, one can consider {\em families} of abelian varieties, and define a 
transform that takes (complexes of) sheaves on a family of abelian varieties 
to (complexes of) sheaves on the corresponding {\em dual} family. 
This has been applied with great success to the study of stable sheaves on 
elliptic surfaces (i.e. holomorphic families of elliptic curves 
parametrized by 
an algebraic curve), see for instance \cite{B} and the references there.
In particular, given an elliptic surface $X$, it can be shown that there 
exists a $1$--$1$ correspondence between vector bundles on $X$ which are 
stable with respect to some suitable polarization, and certain torsion 
sheaves (spectral data) on the relative Jacobian surface 
(see \cite{HRMP, JM} for details).

On the other hand, Mukai's construction can be generalized from complex tori 
to {\em real} tori. Such a generalization, first considered by 
Arinkin \& Polishchuk in \cite{AP} is briefly described in 
Section \ref{locsys} below.
 
Building on previous work by Arinkin \& Polishchuk \cite{AP} and by Bruzzo, 
Marelli \& Pioli \cite{BMP1, BMP2}, we consider in this 
paper a Fourier--Mukai transform for vector bundles with (partial) connections on families of real tori. 
Rather than restricting ourselves to flat connections on {\em Lagrangian} 
families of real tori as in \cite{AP, BMP2},
 we study a broader class of 
connections on vector bundles over a (not necessarily symplectic) manifold 
$M$ with the structure of a {\em flat} torus bundle.

After briefly reviewing the Fourier--Mukai transform for real tori, following 
\cite{BMP1}, we start in Section \ref{fm1} by defining a Fourier--Mukai 
transform for {\em foliated bundles}, which in our context can be viewed as 
{\em families of flat bundles} on the fibers of the torus bundle $M \to B$. 
This takes foliated Hermitian vector bundles over $M$ into certain torsion 
sheaves on the dual fibration $\what{M} \to B$. We then introduce the concept 
of {\em Poincar\'e basic connections} in Section \ref{fm2}, and extend our 
construction to include vector bundles provided with such connections. 
We then conclude in Section \ref{examples} by applying our techniques to 
three different examples: flat connections, that is unitary local systems, 
instantons on $4$--dimensional circle fibrations, and monopoles on 
3-dimensional circle fibrations. 

It is a somewhat surprising fact that certain concepts and techniques 
from foliation theory occur quite naturally in the context of the 
Fourier--Mukai transform. Besides the notions of foliated bundle and 
Poincar\'e basic connection which refer to the torus fibration, that 
is to a foliation which is rather trivial from the point of view 
of foliation theory, there is also a canonical foliation on the dual 
fibration $\what{M} \to B$, transverse to the fibers which has a more 
complicated structure. 
For locally trivial families of flat bundles on the fibers of the torus 
bundle $M \to B$, it turns out that the supports $\Sigma \hra \what{M}$ 
of the Fourier--Mukai transform are (finite unions of) leaves of this 
transverse foliation. This allows us to give in section \ref{flat} an 
explicit parametrization of the representation variety $\cR_M(n)$ of $M$ 
in terms of leaves with transversal holonomy of order $\ell$ such that 
$\ell \vert n~.$ 

The main motivation behind \cite{AP} and \cite{BMP1, BMP2} 
comes from string theory and the Strominger-Yau-Zaslow approach to mirror symmetry, 
with the main goal of understanding Kontsevich's homological mirror symmetry 
conjecture. In a sense, the two main results here presented may also be 
relevant to the understanding of Kontsevich's conjecture. 
Although it seems reasonable to expect that the ideas explored in this 
paper might provide some interesting connections with String Theory 
and mirror symmetry, we do not elaborate on them, leaving such a task 
to mathematical physicists. 

\paragraph{Notation.}
We work on the category of real $C^\infty$-manifolds. By a vector bundle over a manifold $X$, we mean
a  $C^\infty$ vector bundle over $X$. We will also identify a vector bundle with the corresponding sheaf
of  $C^\infty$ sections. By the same token, a sheaf on $X$ should be understood as a sheaf of modules over
the algebra of $C^\infty$ functions on $X$.


\section{Local systems on tori} \label{locsys}

Let us begin by briefly recalling the Fourier--Mukai transform for
real tori, as defined by Arinkin and Polishchuk \cite{AP} and
Bruzzo, Marelli and Pioli \cite{BMP1}; the interested reader
should refer to these papers for the details of this construction.

\medbreak
Let $T$ be the $d$--dimensional real torus, that is
$T=\bR^d/\Lambda$ for the rank $d$ integral lattice
$\Lambda\subset\bR^d$~. The associated dual torus is defined as
$\what T=(\bR^d)^*/(\Lambda)^*$, where
\begin{equation}
(\Lambda)^* := \{ z \in (\bR^d)^* : ~ z (y)  \in \bZ~,~ 
\forall y \in \Lambda \}~,
\end{equation}
is the dual lattice. From the exact sequence
\begin{equation}\label{picext}
0 \to  {\Hom_{\bZ}(\Lambda, \U(1))} \lra H^1(T, \cO^*_{T})
\osetl{c_1} {H^2(T,\bZ}) \to 0  ~,
\end{equation}
we see that, up to gauge equivalence, points in $\what T$
parametrize flat unitary connections on the trivial line bundle
$\underline{\bC} = T \times \bC \lra T$, since we have~:
\begin{equation}\label{rep}
\what{T} = H^1(T, \bR)/{H^1(T, \bZ)} \ocong{\exp} 
{\Hom_{\bZ}(\Lambda, \U(1))} \cong \U(1)^d~.
\end{equation}
For $\xi \in \what{T}~,~ x \in \bR^d~,~ a \in \Lambda$ and
$\lambda \in \bC~,$ consider the equivalence relation
\begin{equation}\label{Pequiv}
\begin{aligned}
\bR^d \times \what{T} \times \bC &\lra  \bR^d \times \what{T}
\times \bC/{\sim}~, \\ 
(x + a, \xi, \lambda) &\sim (x, \xi, \exp (\xi(a)) \lambda)~.
\end{aligned}
\end{equation}
The quotient space under `$\sim$' defines the Poincar\'e line
bundle $P\lra T \times \what{T}$~. Let $p$ and $\hat p$ denote the
natural projections of $T\times \what{T}$ onto its first and
second factors, respectively. In accordance with \eqref{Pequiv},
the bundle $P$ has the property that for $\xi \in \what{T}$, the
restriction $P \mid \hat{p}^{-1}(\xi) \cong L_{\xi}$~, where the
latter denotes the flat line bundle parametrized by $\xi$~. It is
straightforward to see that
\begin{equation}
\Omega^1_{T \times \what{T}} = 
p^*\Omega^1_{T}\oplus\hat{p}^*\Omega^1_{\what{T}}~.
\end{equation}
Corresponding to the definition of $P$ and its above property, it
is shown in \cite{BMP1} that there exists a canonical connection
$\nabla_P: P \lra P \otimes\Omega^1_{T\times \what{T}}$~, whose
connection form is given by
\begin{equation}\label{p-conn}
\bA =  2 \pi \iota ~\sum_{j=1}^d \xi_j dz^j~,
\end{equation}
where $\{z^j\}$ are (flat) coordinates on $T$ and $\{\xi_j\}$ are
dual (flat) coordinates on $\what T$~. The connection $\nabla_P$
splits as the sum $\nabla_P^r\oplus\nabla_P^t$~, where
\begin{equation}\label{connections}
\nabla_P^r = (\id_P \otimes r)\circ \nabla_P ~~,~~\nabla_P^t =
(\id_ P \otimes t)\circ \nabla_P~,
\end{equation}
with natural maps $r:\Omega^1_{T\times \what{T}}\lra
p^*\Omega^1_{T}$, and $t:\Omega^1_{T\times \what{T}}\lra
\hat{p}^*\Omega^1_{\what{T}}$~.

\medbreak
For later purposes we shall denote the dual of any
complex vector bundle $E$ by $E^{\vee}$ and in particular the dual
line bundle of $P$ by $\Pv$~.

\medbreak
Now consider the categories $\sky$ and $\loct$ defined
as follows (see \cite{BMP1})~:
\begin{itemize}
\item $\loct$ is the category of {\it unitary local systems} on $T$~.
Its objects are pairs $(E,\nabla)$ consisting of a smooth complex
vector bundle $E\lra T$ and a flat unitary connection $\nabla$~.
Morphisms are simply bundle maps compatible with the connections.
\item $\sky$ is the category of skyscraper sheaves on
$\what T$ of finite length, that is, $\dim H^0(\what T,S)<\infty$,
for all $S\in \Ob(\sky)$~.
\end{itemize}

The {\it Fourier--Mukai transform} is the invertible functor
\begin{equation}
\bfF:\loct\lra\sky~, \end{equation} which we now describe.
Given $(E,\nabla)\in \Ob (\loct)$, we have the 
{\it relative connection}
\begin{equation}\label{rel1}
\begin{aligned}
\nabla_E^r &: p^*E\ot\Pv \lra p^*E \ot \Pv \ot 
p^*\Omega^1_{T}~, \\ 
\nabla_E^r &= (\id_{(E \ot \Pv)} \ot r) \circ 
( \nabla \ot \id_{\Pv}+\id_{E} \ot\nabla_{\Pv})~,
\end{aligned}
\end{equation}
and {\it the transversal connection}
\begin{equation}\label{trans1}
\begin{aligned}
\nabla_E^t &: p^*E \ot \Pv \lra p^*E \ot \Pv \ot 
p^*\Omega^1_{\what T}~, \\ 
\nabla_E^t &= (\id_{(E \ot \Pv )}\ot t) \circ 
(\nabla \ot \id_{\Pv}+\id_{E} \ot \nabla_{\Pv})~.
\end{aligned}
\end{equation}
As a section of 
$\End_s (p^*E \otimes \Pv) \ot \Omega^2_{T\times \what T}$, 
the commutator satisfies 
(see e.g.~\cite{Marelli})~:
\begin{equation}\label{anticom}
\nabla_E^r \nabla_E^t + \nabla_E^t \nabla_E^r =
\id_E \ot \nabla^2_{\Pv}~.
\end{equation}

\begin{lem} \label{one}{\rm \cite{BMP1}}
If  $(E,\nabla)\in \Ob (\loct)$, then$~:$ 
\begin{itemize}
\item[$(1)$]
$\cR^j\hat{p}_*(\ker\nabla_E^r)=0$~,~{\text{for}}~ 
$0 \leq j \leq d-1$~.

\item[$(2)$]
$S = \cR^d\hat{p}_* (\ker\nabla_E^r) \in \Ob (\sky)$~.
\end{itemize}
Moreover, $\dim H^0(\what T,S)= \rank ~E$~.
\end{lem}
We say that $S=\bfF(E,\nabla)$ is the Fourier--Mukai
transform of the local system $(E,\nabla)$~.

\medbreak
Conversely, take $S \in \Ob (\sky)$, and let $\sigma$
be the support of $S$~. Clearly, $\hat{p}^*S\otimes P$ as a sheaf
on $T\times\what{T}$, is supported on $T\times\sigma$~. Thus
\begin{equation}\label{vanish0}
\cR^j p_*(\hat{p}^*S \ot P)=0~, ~{\text{for}}~ 0 < j \leq d~,
\end{equation}
while $E=\cR^0 p_*(\hat{p}^*S \ot P)$ is a locally--free sheaf 
of rank $\dim H^0(\what T,S)$~. In order to get a connection on $E~,$
consider again the relative connection~:
\begin{equation}
\id_S\otimes\nabla_P^r : \hat{p}^* S \ot P \lra
\hat{p}^*S \ot P \ot p^*\Omega^1_T~.
\end{equation}
Pushing it down to $T$, we get a connection
\begin{equation}
\nabla = \cR^0p_*(\id_S \ot\nabla_P^r) : E \lra
E\otimes\Omega^1_T ~,
\end{equation}
since $\cR^0p_*(\hat{p}^* S \ot P \ot p^*\Omega^1_T) = E \ot \Omega^1_T~,$ 
by the projection formula. Since $(\nabla_P^r)^2 = 0~,$ we conclude 
that $\nabla$ is indeed flat, hence $(E,\nabla)\in \Ob (\loct)$, 
as desired. We use the notation $(E,\nabla) = \what{\bfF}(S)$~.

\medbreak
In summary, referring once more to \cite{BMP1}, we
have~:
\begin{prop} \label{inverse.functors}
The functors $\bfF$ and $\what{\bfF}$ are inverse to
each other, and yield an equivalence between the categories
$\loct$ and $\sky$~.
\end{prop}


\section{The Fourier--Mukai transform} \label{fm1}

Let $M$ be a smooth manifold of dimension $m$, which is the total
space of a $d$--torus bundle over a $(m - d)$--dimensional 
connected manifold $B~,$ that is 
\begin{equation}\label{torus.bundle1}
T^d \hra M \osetl{\pi} B~. 
\end{equation}
Given a point $b\in B~,$ we define $T_b = \pi^{-1}(b)$ to be the fiber 
over $b~,$ where the point $o(b)$ marks the origin of $T_b~.$ 
Regarded as a bundle of groups, $\pi : M \lra B$ admits a discrete 
structure group $\Aut(T) \cong \GL(d, \bZ)~,$ and so the former has 
the structure of a flat fiber bundle and admits a $0$--section 
$o : B \lra M~.$ 
Since the fiber $T$ is compact, this flat structure is 
determined by a holonomy homomorphism 
$\rho : \pi_1 (B) \to \Aut (T) \cong \GL(d, \bZ)$ as a twisted product 
\begin{equation}\label{gen.flat1}
M \cong \wti{B} \times_{\rho} T~.
\end{equation}

\begin{rem}\label{diffremark} 
For the purpose of this paper, we may weaken the structure of the fiber 
bundle $\pi : M \to B$ as follows. 
Let $\Diff (T, o)$ be the group of diffeomorphisms of $T$ which fix 
the origin. Then $\Diff (T)$ is given as a crossed product 
$\Diff (T) \cong T \times_{\vphi} \Diff (T, o)~,$ where $T$ acts by 
translations and $\Diff (T, o)$ acts on $T$ in the obvious way. 
Moreover the canonical homomorphism 
$\Diff (T, o) \to \pi_0 (\Diff (T, o))$ to the mapping class 
group is realized as a deformation retraction 
\begin{equation}\label{torus.bundle2}
\begin{CD}
\Diff (T, o)  @>{\what{}}>> \Aut (\Lambda)  \\
@A{\subset}AA   @A{\cong}AA   \\
\Aut (T)  @>{\cong}>> \GL(d, \bZ) 
\end{CD}
\end{equation}
via $\vphi \bra \wti\vphi \bra \what\vphi = \wti\vphi \vert \Lambda 
\in \Aut (\Lambda)~,$ where $\wti\vphi$ is the unique equivariant 
lift of $\vphi \in \Diff (T, o)$ to $\Diff (\bR^d, o)$ and 
$\what\vphi$ coincides with the automorphism $\vphi_*$ induced by 
$\vphi$ on the fundamental group $\pi_1 (T, o) \cong \Lambda~.$ 
The statement about the deformation retraction follows from the fact 
that any diffeomorphism (actually any homeomorphism) which fixes 
the lattice $\Lambda$ is isotopic to the identity, and in fact 
the connected component $\Diff_e (T, o)$ is contractible to the 
identity; an elementary result which is stated in the 1960's thesis 
of John Franks (as pointed out to us by Keith Burns).
This said, we may start with a fiber bundle $\pi:M \to B$ with 
structure group $\Diff (T, o)~.$ This still guarantees the 
existence of the section $o:B \to M$ and the previous holonomy 
homomorphism $\rho:\pi_1 (B) \to \Aut (T)$ is now recovered as 
the canonical homomorphism $\pi_1 (B) \to \pi_0 (\Diff (T, o))$ 
associated to the fiber bundle $\pi:M \to B~.$ 
Formula \eqref{cross.prod1} for $\pi_1 (M)$ remains valid, as well 
as the flat structure \eqref{gen.flat2} of the dual fiber bundle 
$\hat\pi:\what{M} \to B~,$ the latter property being a consequence 
of the homotopy invariance of singular cohomology. 
In fact, the above deformation retraction implies that the structure 
group of a $\Diff (T, o)$--torus bundle admits a unique reduction to 
$\Aut (T) \subset \Diff (T, o)$ and so $\pi:M \to B$ is still 
diffeomorphic to a flat fiber bundle of the form \eqref{gen.flat1}. 
\end{rem}

The fundamental group of $M$ is determined as a crossed product
\begin{equation}\label{cross.prod1}
0 \to \Lambda = \pi_1 (T) \lra \pi_1 (M) 
\cong \pi_1 (T)\times_{\rho_*} \pi_1 (B) \osetl{\la} \pi_1 (B) \to 1~,
\end{equation}
where $\rho_*$ is given by the induced action of $\pi_1 (B)$ on 
$\Lambda$ via the isomorphism $\GL(d, \bZ) \cong \Aut (\Lambda)~.$ 

We have the exact sequence
\begin{equation}\label{exact1}
0 \to T(\pi) \lra TM \lra \pi^* TB \to 0~,
\end{equation}
and the dual sequence of $1$--forms
\begin{equation}\label{exact2}
0 \to \pi^* \Omega^1_B \lra \Omega^1_M \lra \Omega^1_{M/B} \to 0~.
\end{equation}
Observe that a flat structure of $\pi : M \to B$ defines a 
splitting of the exact sequences \eqref{exact1} and \eqref{exact2}. 

The {\em dual fiber bundle} $\what{M} \to B$ is given by the 
total space of $\cR^1\pi_*\bR/ \cR^1\pi_*\bZ$ as a 
(locally constant) sheaf on B~. 
If $\hat{\pi} : \what{M}\lra B$ is the natural projection, it 
is easy to see that $\hat{\pi}^{-1}(b)=\what{T}_b$~. Note that 
this projection also has a $0$--section $\sigma_0:B \to \what{M}~.$ 
It follows that $\hat\pi : \what{M} \to B$ is given by the flat 
bundle of fiber cohomologies 
\begin{equation}\label{gen.flat2}
\what{M} \cong \wti{B} \times_{\rho^*} \what{T}~, 
\end{equation}
where $\rho^*$ is the induced action of $\pi_1 (B)$ on 
$\what{T} = H^1(T, \bR)/{H^1(T, \bZ)}~.$ 
Furthermore, $\cR^1\hat{\pi}_*\bR/ \cR^1\hat{\pi}_*\bZ$ coincides with
$M$ as sheaves on $B$, and we have $\what{\what{M}} \cong M$~.

\medbreak
Let $Z = M \times_B \what M$ be the fiber product, with
its natural projections $p:Z\lra M$ and $\hat p:Z\lra\what M$ onto
the first and second factors. Clearly, 
$\pi \circ p = \what\pi \circ \hat p$ and 
$(\pi\circ p)^{-1}(b)=T_b\times\what{T}_b~.$
\begin{equation}\label{diamond1}
\xymatrix{ & Z \ar[dl]^p \ar[dr]^{\hat p} & \\ 
M \ar[dr]^\pi & & \what{M} \ar[dl]^{\hat\pi} \\ 
& B & }
\end{equation}
It is also easy to see that $p^{-1}(x)=\what{T}_{\pi(x)}$, for all
$x\in M$ and $\hat{p}^{-1}(y)= T_{\hat{\pi}(y)}$, for all
$y\in\what M$~.

\medbreak
Defining $\Omega^1_{Z/\what M}=\Omega^1_Z/\hat{p}^*\Omega^1_{\what M}~,$ 
recall that the Gauss--Manin connection yields a splitting of the short exact sequence 
\begin{equation}\label{exact3}
0 \to \hat{p}^*\Omega^1_{\what M} \lra \Omega^1_Z
\stackrel{r}{\lra} \Omega_{Z/\what M}^1 \to 0~, 
\end{equation}
such that we have the decomposition
\begin{equation}\label{split1}
\Omega^1_Z = \hat p^*\Omega^1_{\what M} \oplus \Omega^1_{Z/\what{M}}~.
\end{equation}
From (\ref{exact2}) it follows that 
\begin{equation}\label{fibremark}
\Omega^1_{Z/\what{M}} = p^* \Omega^1_{M/B}~,
\end{equation}
since $\hat{p}: Z \to \what{M}$ is the pull--back 
fibration of $\pi : M \to B$ along $\hat{\pi}~.$

\medbreak
There exists a line bundle $\rpb$ over $Z = M \times_B \what M$,
with the property that $\rpb \vert {(\pi\circ p)^{-1}(b)}$ is just
the Poincar\'e line bundle $\cP_b$ over $T_b\times\what T_b$, for
all $b\in B$ (see \cite{BMP2}). We call $\rpb$ {\it the relative
Poincar\'e line bundle}. Just as in the absolute case, it has the
property that for $\xi \in \what{M}$, the restriction $\rpb \mid
\hat{p}^{-1}(\xi) \cong L_{\xi}$, where the latter denotes the
flat line bundle parametrized by $\xi \in \what{T}_b \subset
\what{M}~.$ 

There is a canonical connection on $\rpb$ which we
denote by $\nabla_\rpb$~. Following \cite{BMP2}, we can write 
its connection matrix $\bA$ in a suitable gauge on an open
subset $U \times T \times \what{T} \subset Z$ 
as follows
\begin{equation}\label{rp-conn}
\bA =  2 \pi \iota ~\sum_{j=1}^d \xi_j dz^j~,
\end{equation}
where $\{z^j\}$ are (flat) coordinates on $T$ and $\{\xi_j\}$
are dual (flat) coordinates on $\what T$~. In such coordinates
the curvature
$\bF = \nabla_{\rpb}^2$ is then given by~:
\begin{equation} \label{p-curv}
\bF = 2\pi \iota~ \sum_{j=1}^d  d\xi_j \wedge dz^j~.
\end{equation}
In the same coordinate system, we have
$\nabla_{\rpbv}^2 = -\bF$~.

\subsection{Transforming foliated bundles}
\label{foliated}

Let $E \ra M$ be a Hermitian vector bundle of rank $n$~. With
reference to \eqref{exact1}, we assume a {\em foliated bundle} 
structure on $E$ given by a {\em flat partial unitary connection} 
\cite{KTone}~:
\begin{equation}
\npart_E : E \lra E \ot \Omega^1_{M/B} = E \ot \Omega^1_M/
\pi^* \Omega^1_B~,
\end{equation}
satisfying $(\npart_E)^2 = 0$~.

\medbreak
The local structure of a foliated bundle $(E, \npart_E)$ on $M$ 
is described next. 
 
\begin{ex}\label{holonomy1}
Local structure of foliated bundles on $M~:$~ 
Intuitively, a foliated bundle on $M \oset{\pi} B$ is a family 
of flat bundles (unitary local systems) $(E_b,\nabla_{E_b})$ on 
the fibers $T_b~,$ parametrized by $b \in B~.$ Of course, the 
topology of $E$ has to be taken into account. The local description 
is quite similar to that of the Poincar\'e line bundle in 
\eqref{Pequiv}, which is of course an example of a foliated bundle. 
Thus for sufficiently small open sets $U \subset B~,$ there are 
isomorphisms
\begin{equation}\label{holonomy2}
\begin{CD}
U \times (\bR^d \times_{\Lambda} \bC^n) @>{\cong}>> 
E \mid \pi^{-1} (U) \\
@VV{\Id \times \tau'}V    @VV{\tau\vert\pi^{-1}U}V  \\
U \times T   @>{\cong}>>   \pi^{-1} (U)~, 
\end{CD}
\end{equation}
where the identification on the LHS is given by 
\begin{equation}\label{holonomy3}
(b, x + a, \lambda) \sim (b, x, \exp (\xi_b (a)) \lambda)~, 
\end{equation}
for $\xi = (\xi_1, \ldots, \xi_n)~, ~\xi_j : U \to \what{T}~, ~ 
\what{T} \ocong{\exp} \Hom_{\bZ}(\Lambda, \U(1))~, ~b \in U~, ~
x \in \bR^d~,~ a \in \Lambda$ and $\lambda \in \bC^n~.$ 
Relative to a (good) open cover $\cU$ of $B~,$ we have coordinate 
changes over $U_{ik} = U_i \cap U_k$ of the form 
$$
(\Id, \wti\vphi_{ik}, g_{ik}) : 
U_{ik} \times \bR^d \times \bC^n \osetl{\cong} 
U_{ik} \times \bR^d \times \bC^n~, 
$$
compatible with the identifications in \eqref{holonomy3}, that is 
\begin{equation}\label{holonomy4}
\Adj (g_{ik} (b)) \circ \exp(\xi^k_b (a)) = 
\exp(\xi^i_b (\what\vphi_{ik}~a))~. 
\end{equation}
Here $\{ \vphi_{ik} \}$ is the smooth $1$--cocycle on $\cU$ with 
values in $\Diff (T, o)$ describing the fiber bundle $\pi:M \to B~, ~
\wti\vphi$ is the unique equivariant lift of 
$\vphi \in \Diff (T, o)$ to $\Diff (\bR^d, o)$ and 
$\what\vphi = \wti\vphi \vert \Lambda$ 
is the induced automorphism on the lattice $\Lambda~.$ 
$\{ g_{ik} \}$ is a smooth $1$--cochain of local gauge 
transformations $g_{ik} : U_{ik} \to \U (n)$ on $\cU~.$ 
\end{ex}

\medbreak
From \eqref{p-curv} we see that the unitary connections 
$\nabla_{\rpb}$ and $\nabla_{\rpbv }$ are flat along the 
fibers of the projection $\hat{p} : Z \to \what{M}$ in 
\eqref{diamond1} and induce flat partial unitary connections 
$\npart_{\rpb}$ and $\npart_{\rpbv}$ on $\rpb$ and $\rpbv$~.
Pulling $E$ back to $Z$ and tensoring with the dual Poincar\'e
bundle $\rpbv$~, consider the flat partial connection~:
\begin{equation} \label{tilde.conn1}
\wti{\nabla}^r_E = p^* \npart_E \ot \id_{\rpbv} +
\id_E \ot \npart_{\rpbv} : p^*E \ot \rpbv \lra
p^* E \ot \rpbv \ot \Omega^1_{Z / \what{M}}~.
\end{equation}
Now for each $b \in B$, the pair $(E,\npart_E)$ restricts
to a unitary local system $(E_b,\nabla_{E_b})$ over the fiber
$T_b$, while the connection $\wti{\nabla}^r_E$ restricts to the
operator $\wti{\nabla}_{E_b}^r$ induced by \eqref{rel1}.
Therefore,
\begin{equation} \label{iso1}
\cR^j \hat{p}_{b,*} \left( (\ker \wti{\nabla}^r_E)
\vert{T_b\times\what{T}_b} \right) \cong \cR^j \hat{p}_{b,*} 
( \ker \wti{\nabla}_{E_b}^r )~,
\end{equation}
where $\hat{p}_{b}:T_b\times\what{T}_b\lra\what{T}_b$ is the
projection onto the second factor.

\medbreak
On the other hand, let $\iota_b,\hat{\iota}_b$ be the
inclusions of $T_b~, \what{T}_b$ into $M$ and $\what{M}$
respectively, and consider the diagram~:
\begin{equation}
\xymatrix{ T_b\times\what{T}_b \ar[r]^{\iota_b\times\hat{\iota}_b}
\ar[d]^{\hat{p}_b} & Z \ar[d]^{\hat{p}} \\ \what{T}_b
\ar[r]^{\hat{\iota}_b} & \what{M} }
\end{equation}
Then the topological base change \cite{Gelf} yields the
isomorphism (for $0 \leq j \leq d$)~:
\begin{equation}
\left. \cR^j \hat{p}_{*} (\ker \wti{\nabla}^r_E) \right
\vert {\what{T}_b} \cong \cR^j \hat{p}_{b,*} \left( (\ker
\wti{\nabla}^r_E) \vert {T_b\times\what{T}_b} \right)~.
\end{equation}
Combining with (\ref{iso1}), one obtains~:
\begin{equation} \label{iso2}
\left. \cR^j \hat{p}_{*} (\ker \wti{\nabla}^r_E) \right \vert
{\what{T}_b} \cong \cR^j \hat{p}_{b,*} (\ker \wti{\nabla}_{E_b}^r)~.
\end{equation}
It then follows from Lemma \ref{one} that 
\begin{equation}\label{vanish1}
\cR^j \hat{p}_{*} (\ker \wti{\nabla}^r_E) = 0~, ~j=0,\ldots,d-1~. 
\end{equation}
Now we set
\begin{equation}\label{transf.bdl1}
\what{E} = \cR^d \hat{p}_{*} (\ker \wti{\nabla}^r_E)~,
\end{equation}
with its support denoted by 
$\Sigma = \Sigma (E, \npart_E) = \supp \what{E}$~. 
From (\ref{iso2}), we have then
\begin{equation} \label{iso3}
\what{E} \mid {\what{T}_b} \cong \cR^d \hat{p}_{b,*} 
( \ker \wti{\nabla}_{E_b}^r )~.
\end{equation}

\medbreak
The following elementary Lemma is useful. 
As it is purely local, it is valid for any foliated bundle 
$E \to M~.$

\begin{lem}\label{parallel.sec}
For sufficiently small open sets $V \subset M~,$ the foliated 
Hermitian vector bundle $E \mid V$ admits $\npart_E$--parallel 
unitary frames $s = (s_1, \ldots, s_n)~,$ that is 
$\npart_E ~s_i = 0~, ~i=1,\ldots,n~.$ 
It follows that $\npart_E$ is linear over the sheaf $\pi^* \cO_B$ 
of basic functions and the sheaf $\ker \npart_E$ of 
$\npart_E$--parallel sections is locally free as a module over 
$\pi^* \cO_B~,$ of the same rank as $E~.$ 
\end{lem}

\begin{proof}
This can be shown easily by working in a sufficiently small 
Frobenius chart $V = U \times U' \subset U \times T$ over 
which $E$ trivializes, choosing any unitary frame along $U \times a~, 
~a \in U'$ and then using parallel transport relative to the flat 
partial unitary connection $\npart_E$ in the fiber direction $U'~.$ 
\end{proof}
In our context, this means that $\wti{\nabla}^r_E$ is linear with
respect to $\hat{p}^* \cO_{\what M}$~. In particular, the sheaf 
$\ker \wti{\nabla}^r_E$ of $\wti{\nabla}^r_E$--parallel sections in 
$p^*E \otimes \rpbv$ is a locally free module over 
$\hat{p}^* \cO_{\what M}~.$ 
Further, the derived direct image $\what E$ is a torsion module 
over $\cO_{\what M}$~.

\medbreak
For $\Sigma = \supp \what{E}~,$ we shall also 
consider the fiber product $Z_\Sigma = M\times_B \Sigma$, 
with $p_\Sigma: Z_\Sigma \lra M$, and 
$\hat{p}_\Sigma: Z_\Sigma\lra \Sigma$, 
denoting the natural projections~:
\begin{equation}
\xymatrix{ & Z_{\Sigma} \ar[dl]^{p_{\Sigma}} \ar[dr]^{\hps} & \\ 
M \ar[dr]^\pi & & \Sigma \ar[dl]^{\hpis} \\ 
& B & }
\end{equation}
There is also the restriction to $Z_{\Sigma}$ of the relative
Poincar\'e line bundle $\rpb$; we denote this by $\rpb_\Sigma =
\rpb \vert Z_{\Sigma}$~. Let $j:\Sigma \hra \what{M}$ be
the inclusion map, and let $\tilde{j}:Z_\Sigma \hra Z$
be the induced inclusion. 

\medbreak
Next we set $\cK = \tilde{j}^* (\ker \wti{\nabla}^r_E)$,
and consider the sheaf
\begin{equation}\label{supp1}
\cL = j^* \what{E} = \cR^d\hpsstar(\cK)~.
\end{equation}

\begin{prop} \label{two}
For $\what{E}$ given by \eqref{transf.bdl1} and 
$\Sigma = \supp \what{E}~,$ we have$~:$

\begin{itemize}
\item[$(1)$]
$\what{E} \mid \what{T}_b \in \Ob (\skynb)$ for $b \in B$ and 
the support $\Sigma$ of $\what{E}$ is closed and transversal to 
all fibers $\what{T}_b~.$ 

\item[$(2)$]
For $\Sigma_b = \Sigma\cap\what{T}_b = \supp (\what{E} \mid \what{T}_b)~,$ 
the counting function $\abs{\Sigma_b}$ satisfies 
$1 \leq \abs{\Sigma_b} \leq n~, \forall b\in B~.$ 
The sets $U_\ell \subset B~, ~\ell = 1, \ldots, n~,$ 
for which $\abs{\Sigma_b} \geq \ell$ are open in $B,$ 
possibly empty for $\ell > 1~,$ and satisfy 
$U_n \subseteq \ldots \subseteq U_{\ell+1} 
\subseteq U_{\ell} \subseteq \ldots \subseteq U_1 = B~.$  

\item[$(3)$]
For every $b \in B~,$ there is an open neighborhood 
$U \subset B$ of $b~,$ such that the con\-nec\-ted 
com\-po\-nents $\hpis^{-1} (U)_{\xi}$ of 
$\hpis^{-1} (U) \subset \Sigma$ containing $\xi \in \Sigma_b$ 
separate the elements $\xi \in \Sigma_b$ and $\hpis^{-1} (U)_{\xi}$ 
can be exhausted by a finite number of smooth sections 
$\sigma_i : U \to \hpis^{-1} (U)_{\xi}~,$ such that 
$\sigma_i (b) = \xi~.$ For $U$ sufficiently small, 
the number of sections needed is bounded by the rank 
of $\cL$ at $\xi \in \Sigma_b~.$ 

\item[$(4)$]
The rank of $\cL \lra \Sigma$ at $\xi \in \Sigma_b$ is
equal to the multiplicity of the irreducible representation 
$\exp(\xi)$ in the unitary local system $(E_b,\nabla_{E_b})$ 
on $T_b~,$ that is the multiplicity of the trivial representation
in the flat bundle $E_b \otimes L_\xi^{\vee} \lra T_b~.$
\end{itemize}
\end{prop}

We say that $\Sigma \hra \what{M}~,$ satisfying $(1)$ to $(3)$ 
in Lemma \ref{two}, is a {\em $n$--fold ramified covering of $B$ 
of dimension $m-d = \dim (B)$}. 
A point $\xi \in \Sigma$ is called {\em regular} or {\em smooth}, 
if the connected component $\hpis^{-1} (U)_{\xi}$ is given by a 
single section $\sigma:U \cong \hpis^{-1} (U)_{\xi}$ for a 
sufficiently small open neighborhood $U$ of $b = \hat\pi (\xi)~.$ 
The {\em regular} set $\Sigma_{\reg} \subseteq \Sigma$ is the set 
of regular points in $\Sigma~.$ $\Sigma_{\reg}$ is an open, dense 
subset of $\Sigma~,$ $\Sigma_{\reg} \hra \what{M}$ is a smooth 
submanifold and the rank of $\cL$ is locally constant on 
$\Sigma_{\reg}~,$ that is $\cL$ is a locally free module on the 
connected components of $\Sigma_{\reg}~.$ 
The closed, residual complement 
$\Sigma_{\br} = \Sigma \bsl \Sigma_{\reg}$ 
is called the {\em branch locus} of $\hpis : \Sigma \to B~.$ 
We say that $\Sigma \hra \what{M}$ is {\em smooth} if the 
branch locus $\Sigma_{\br}$ is empty. In this case we have 
$\Sigma_{\reg} = \Sigma$ and $\Sigma \hra \what{M}$ is a closed 
smooth submanifold, the rank of $\cL$ is locally constant on $\Sigma$ 
and the semicontinuous counting function $\abs{\Sigma_b}$ is locally 
constant, hence constant on $B~.$ 
In particular, $\Sigma$ is smooth if $U_n = B~,$ that is 
$\abs{\Sigma_b} \equiv n$ on $B~,$ in which case $\cL$ is a complex 
line bundle on $\Sigma~.$ If $\Sigma$ is in addition connected, 
then $\Sigma$ is a smooth $n$--fold covering space of $B$ in the 
usual sense and we say that $\hpis :\Sigma \to B$ is 
{\em non--degenerate}. 

\begin{proof}
Lemma \ref{one} and the identification \eqref{iso2} imply that
\begin{equation}\label{transf.bdl2}
\what{E} \mid \what{T}_b  = \cR^d\hat{p}_{b,*} 
( \ker \wti{\nabla}_{E_b}^r ) \in \Ob (\skynb)~.
\end{equation}
Since $\Sigma_b = \supp (\what{E} \mid \what{T}_b)~,$ and 
$\dim H^0(\what{T}_b, \what{E} \mid \what{T}_b)= n~,$ 
part $(1)$ follows easily.

For $\xi \in \Sigma_b$~, we have $\hat{p}^{-1} (\xi) = T_b$ and 
from \eqref{iso3}, we see that the rank of $\cL$ at $\xi$
is given by the rank of the cohomology group 
$H^d (T_b, \ker \wti{\nabla}^r_{E_b})$~. This proves $(4)$. 

From $(1),$ we have $1 \leq \abs{\Sigma_b} \leq n~.$ 
From $(4),$ we see that the second condition in $(2)$ is really 
the semicontinuity of the number of distinct holonomy 
representations $\xi \in \Sigma_b$ in the bundles $E_b~.$ 
Thus for $b \in B~,$ there is a neighborhood $U_b \subset B~,$ such that 
$\abs{\Sigma_{b'}} \geq \abs{\Sigma_b}~, ~b' \in U_b~.$ Then 
$b \in U_{\ell}$ implies that $U_b \subset U_{\ell}$ and $(2)$ follows. 

Finally, $(3)$ is proved by using the local description 
of a foliated bundle in Example \ref{holonomy1} and $(4)~.$ 
In fact, the number of sections needed is equal to the number 
of distinct germs at $b$ among the functions $\xi_j$ passing 
through $\xi$ in \eqref{holonomy3} and therefore is bounded 
by the rank of $\cL$ at $\xi \in \Sigma_b~.$ 
\end{proof}

\medbreak
In view of the above result, we say that $\what{E}$ is a 
{\em relative skyscraper} (that is, a sheaf whose restriction 
$\what{E} \mid \what{T}_b $ to each fiber $T_b$ is a skyscraper 
sheaf of constant finite length), $\cL$ is the 
{\em sheaf of multiplicities} and $\Sigma$ is the 
{\em spectral covering} of $(E, \npart_E)$~. Note that these structures
are completely determined by the flat partial connection $\npart_E$~.

\subsection{The inverse transform for relative skyscrapers} 
\label{inverse1}

The inverse construction is considerably simpler. Our starting 
point is the pair $(\cS, \Sigma)~,$ where $\cS$ is a relative 
skyscraper of constant length $n$ on $\what M$ supported on a 
$n$--fold ramified covering $\Sigma \hra \what{M}$ of $B$ of 
dimension $m-d = \dim (B)~.$ 

\medbreak
Using the same notation as before, recall that the
fiber product $Z_{\Sigma} = M \times_B \Sigma$, is of dimension
$m$, and $p_\Sigma:Z_\Sigma\lra M$ is an $n$--fold ramified
covering map. Thus it is easy to see that
\begin{equation}\label{checks}
\check{\cS} = \psstar (\hps^* \cS \otimes \rpb_\Sigma)~,
\end{equation}
is a locally--free sheaf of rank $n$ on $M$~. Furthermore, the
construction reveals that $\check{\cS}$ carries a canonical 
flat partial connection $\npart_{\check{\cS}}$~. In fact, 
$\hat{p}_{\Sigma}^* \cS$ carries a canonical flat partial 
connection relative to $\hat{p}_{\Sigma} : Z_\Sigma \to \what{M}$ 
and so does $\rpb_\Sigma$~.

\subsection{The main result} \label{main1}

Motivated by the results above, let us introduce the following
categories of sheaves with connections on $M$ and $\what M$~.
\begin{dfn}
$\vectno$ is the category of foliated Hermitian vector bundles 
on $M$ endowed with a flat partial unitary connection. Objects in 
$\vectno$ are pairs $(E,\npart_E)$ consisting of a Hermitian vector 
bundle $E$ of rank $n$ and a flat partial unitary connection $\npart_E$~.
Morphisms are bundle maps compatible with such connections.
\end{dfn}

\begin{dfn}
$\relskyn$ is the category of {\it relative skyscrapers} on 
$\what M$. Objects in $\relskyn$ are pairs $(\cS, \Sigma)$ 
consisting of a relative skyscraper $\cS$ of constant length 
$n$ on $\what M~,$ supported on a $n$--fold ramified covering 
$\Sigma \hra \what{M}$ of $B$ of dimension $m-d = \dim (B)~.$ 
Morphisms are sheaf maps of $\cO_{\what{M}}$--modules. 
\end{dfn}

\medbreak
The constructions in Sections \ref{foliated} and \ref{inverse1} 
define additive covariant functors 
\begin{equation}\label{functors} 
\fm : \vectno \lra \relskyn~,~~ \fmh : \relskyn \lra \vectno~.
\end{equation}
For limits in the appropriate sense, let
\begin{equation}\label{fcats}
\vecto = \lim_n \vectno~~{\text{and}}~~ \relsky = \lim_n
\relskyn~.
\end{equation} 
With these definitions in place, we can state our main result~:

\begin{thm} \label{mainresult1}
The Fourier--Mukai trans\-form $\fm$ defines an additive 
natural equi\-valence of cate\-gories
\begin{equation}
\fm ~:~\vecto \osetl{\cong} \relsky~.
\end{equation}
\end{thm}
\begin{proof}
We claim that $\fm$ and $\fmh$ are
adjoint functors which in fact define an equivalence of
categories. From the construction of $\fm$ and $\fmh$, 
there exist natural transformations
\begin{equation}\label{nat.transf1}
\phi_{\cS} : \cS \lra \fm \circ \fmh (\cS)~, ~ 
\cS \in \Ob(\relskyn)~,
\end{equation}
and
\begin{equation}\label{nat.transf2}
\psi_E : \fmh \circ \fm (E) \lra E ~, ~
E \in \Ob (\vectno)~.
\end{equation}
These natural transformations define adjunction maps
\begin{equation}\label{adjunct}
\begin{aligned}
\Phi  &: \Morph_{\bfV}(\fmh(\cS), E) \lra \Morph_{\bfR}(\cS,
\fm(E))~,\\ 
\Psi &: \Morph_{\bfR}(\cS, \fm(E)) \lra 
\Morph_{\bfV}(\fmh(\cS), E)~,
\end{aligned}
\end{equation}
where $\Morph_{\bfV}$ and  $\Morph_{\bfR}$ denote morphisms in
$\vectno$ and $\relskyn$, respectively. 
Explicitly, for $f : \fmh(\cS) \lra E$~, we have by naturality 
\begin{equation} \label{Phi}
\Phi(f) = 
\fm(f) \circ \phi_{\cS}~,
\end{equation}
so that $\phi_{\cS}$ determines $\Phi$~. 
Likewise, for $g : \cS \lra \fm(E)$, we have by naturality 
\begin{equation} \label{Psi}
\Psi (g) = 
\psi_E \circ \fmh (g)~,
\end{equation}
so that $\psi_E$ determines $\Psi$ as well. 
The natural transformations \eqref{nat.transf1}, \eqref{nat.transf2} 
correspond then to $\phi_{\cS} = \Phi(\id_{\fmh(\cS)})$ and 
$\psi_E = \Psi(\id_{\fm(E)})$, respectively. 
The fact that the adjunction maps
$\Phi$ and $\Psi$ are inverses of each other, is equivalent to the
compositions
\begin{equation}
\begin{aligned}
\fm (E) \osetl{\phi_{\fm(E)}} \fm \circ \fmh \circ (\fm(E)) &= \fm
\circ (\fmh \circ \fm(E)) \osetl{\fm(\psi_E)} \fm(E)~, \\ \fmh
(\cS) \osetl{\fmh(\phi_{\cS})} \fmh \circ (\fm \circ \fmh(\cS)) &=
\fmh \circ \fm \circ( \fmh (\cS)) \osetl{\psi_{\fmh(\cS)}}
\fmh(\cS)~,
\end{aligned}
\end{equation}
resulting in the identities of $\fm(E)$ and $\fmh(\cS)$,
respectively.

\medbreak
The construction has the further property that it is
compatible with localization relative to open subsets $U \subset
B$, that is, the restrictions to  $\pi^{-1}(U)$ and $\hat
\pi^{-1}(U)$~. Moreover, we observe that the restriction of $\fm$
and $\fmh$ to the fibers of $M$ and $\what M$ at $b \in B$ 
respectively, coincides with the functors
\begin{equation}
\fm_b : \loctb \lra \skyb ~~,~~
\what{\fm}_b: \skyb \lra \loctb~,
\end{equation}
for each $b \in B$~. It follows from \cite{BMP1,Marelli}
that $\phi_{\cS_b}:\id_{{\what T}_b} \cong \fm_b \circ \fmh_b$ and 
$\psi_{E_b}:\fmh_b \circ \fm_b \cong \id_{T_b}$~. From this we 
conclude that $\phi_{\cS}$ and $\psi_E$  are indeed isomorphisms.
\end{proof}

\medbreak
Let $V \in \Vct{n}{B}~,$ where $\Vct{n}{B}$ is the category of 
complex vector bundles of rank $n$ over $B~.$ Then $\pi^* V$ 
carries a canonical flat partial connection 
$\npart_{\pi^* V}~,$ so that $(\pi^* V, \npart_{\pi^* V})$ 
is an object in $\Vecto{n}{M}~,$ while 
$\hat{\pi}_0^* V = \hat{\pi}_{\Sigma_0}^* V$ 
is an object in $\relskyn~,$ supported on the 
$0$--section $\Sigma_0 = \sigma_0 (B) \subset \what{M}~.$ 
The construction of $\fm$ is compatible with these pull--backs, 
that is we have a commutative diagram  
\begin{equation}\label{diagram3}
\begin{CD}
\Vecto{n}{M}  @>{\fm}>>  \relskyn   \\
@A{\pi^*}AA   @A{\hat\pi_0^*}AA   \\
\Vct{n}{B}  @>{=}>> \Vct{n}{B} . 
\end{CD}
\end{equation}
Moreover, the Fourier--Mukai transform $\fm$ has a module property  
with respect to $\Vect{}{B}~.$ 

\begin{cor}\label{module1}
For $(E, \npart_E) \in \Vecto{}{M}$ and $V \in \Vct{}{B}~,$ 
the Fourier--Mukai trans\-form $\fm$ satisfies 
\begin{equation}\label{module2}
\fm ((\pi^* V, \npart_{\pi^* V}) \ot (E,\npart_E)) \cong 
\hpis^* V \ot \fm (E, \npart_E)~,
\end{equation}
where $\Sigma$ is the support of $\fm (E,\npart_E)~.$ 
\end{cor}


\section{The Fourier--Mukai transform for vector bundles with  \\ 
Poincar\'e basic connections} \label{fm2}

\subsection{Transforming bundles with Poincar\'e basic connections} 
\label{basic1}

Let $E \lra M$ be a foliated Hermitian vector bundle of rank $n$, 
and let $\nabla_E : E \lra E \otimes \Omega^1_M$ be a unitary 
connection on $E$~. We say that $\nabla_E$ is {\em adapted} 
to the foliated structure on $E$~, if $\nabla_E$ induces the 
flat partial connection $\npart_E : E \lra E \otimes \Omega^1_{M/B}$ 
via the canonical map $\Omega^1_M \to \Omega^1_{M/B}$ in \eqref{exact2}. 
The existence of adapted connections follows from an elementary 
partition of unity argument. 

\medbreak
At this point, it is also useful to introduce the bigrading on 
the DeRham algebra $\Omega^*_M$ determined by a 
splitting of the exact sequence \eqref{exact1}, 
respectively \eqref{exact2}~: 
\begin{equation}\label{bigrad1}
\Omega^{u,v}_M = \Omega^{u,0}_M \otimes 
\Omega^{0,v}_M  = \pi^* \Omega^u_B \otimes \Omega^v_{M/B}~.
\end{equation}
$u$ is called the {\em transversal or basic} degree and $v$ 
is called the {\em fiber} degree. 

\medbreak
Consider now the adapted connection 
\begin{equation} \label{tilde.conn}
\wti{\nabla}_E = p^*\nabla_E \ot \id_{\rpbv} +
\id_E \ot \nabla_{\rpbv} : 
p^*E \ot \rpbv \lra p^*E \ot \rpbv \ot \Omega^1_{Z}~, 
\end{equation}
on $p^* E \ot \rpbv$~. 
With respect to a corresponding splitting of \eqref{exact3}, 
we have $\wti{\nabla}_E =\nabla_E^r \oplus \nabla_E^t$, where
\begin{equation} \label{rel.conn}
\nabla_E^r = (\id_{E\otimes\rpbv}\otimes r)\circ \wti{\nabla}_E :
p^*E\otimes\rpbv \lra p^*E\otimes\rpbv\otimes\Omega^1_{Z/{\what M}}~,
\end{equation}
is the relative connection, and
\begin{equation} \label{trans.conn}
\nabla_E^t = (\id_{E\otimes \rpbv} \otimes t)\circ \wti{\nabla}_E
: p^*E \otimes \rpbv \lra p^*E \otimes \rpbv
\otimes\hat{p}^*\Omega^1_{\what M}~,
\end{equation}
is the transversal connection, that is the components of type 
$(0,1)$ and $(1,0)$ of $\nabla_E$ respectively. 

In the sequel, we always view the curvature $\nabla_E^2$ of 
$\nabla_E$ as a $2$--form with values in the adjoint bundle 
$\End_s (E)$ of skew--hermitian endomorphisms of $E$~.

\medbreak
\begin{lem}\label{anticomlemma}
The type--decomposition of the curvature $\wti{\nabla}_E^2$ is 
given by
\begin{equation}\label{curv0}
(\wti{\nabla}_E^2)^{0,2} = p^*((\nabla_E^2)^{0,2}) 
\ot \id_{\rpbv} = 0~,
\end{equation}
\begin{equation}\label{curv1}
(\wti{\nabla}_E^2)^{2,0} = p^*((\nabla_E^2)^{2,0}) 
\ot \id_{\rpbv} = (\nabla_E^t)^2~,
\end{equation}
and 
\begin{equation}\label{curv2}
(\wti{\nabla}_E^2)^{1,1} = p^*((\nabla_E^2)^{1,1}) \ot 
\id_{\rpbv} - \id_E \otimes \bF = \Xi~, 
\end{equation}
where the operator $\Xi$ is given by the 
commutator 
\begin{equation}\label{anticomm1}
\Xi = \nabla_E^r \circ  \nabla_E^t - \nabla_E^t \circ  \nabla_E^r~:~
p^*E\ot\rpbv \lra p^*E\ot\rpbv\ot\Omega^{1,1}_Z~.
\end{equation}
\end{lem}
Henceforth we adopt the usual sign rule which equips the extension 
of the transversal operator $\nabla_E^t$ to forms of higher degree 
with a sign $(-1)^v$ on forms of type $(u,v)~.$ 

\begin{proof}
Firstly from \eqref{tilde.conn} and the decomposition 
\eqref{rel.conn}, \eqref{trans.conn} we have
\begin{equation}
\wti{\nabla}_E = p^*\nabla_E\otimes\id_{\rpbv} +
\id_E\otimes\nabla_{\rpbv} = \nabla_E^r + \nabla_E^t~.
\end{equation}
Computing the curvature operator $\wti{\nabla}_E^2$ in two ways, 
we obtain~:
\begin{equation}\label{curv.op1}
\begin{aligned}
\wti{\nabla}_E^2 &= (p^*\nabla_E)^2 \otimes \id_{\rpbv} + \id_E
\otimes \nabla_{\rpbv}^2 \\ &= p^*(\nabla_E^2) \otimes \id_{\rpbv} - 
\id_E \otimes \bF \\ &= 
(\nabla_E^r \pm \nabla_E^t)^2   \\ 
&= (\nabla_E^r)^2 + (\nabla_E^t)^2 + \Xi~.
\end{aligned}
\end{equation}
Since $\nabla_E$ is adapted to the foliated structure $\npart_E$ 
on $E$~, we have 
$({\nabla}_E^2)^{0,2} = 0$~. 
Since $\wti{\nabla}_E$ is adapted to the foliated structure  
on $p^*E \otimes \rpbv$ relative to $\hat{p} : Z \to \what{M}$~, 
we have $(\wti{\nabla}_E^2)^{0,2} = (\nabla_E^r)^2 = 0$~. 
The curvature $\bF$ of the relative Poincar\'e bundle $\rpb$ 
is of type $(1,1)$ by \eqref{p-curv} and $(\nabla_E^t)^2$ 
and $\Xi$ are of type $(2,0)$ and $(1,1)$ respectively by definition. 
Thus the assertions \eqref{curv0}, \eqref{curv1} and \eqref{curv2} 
follow from \eqref{curv.op1}. 
We use \eqref{fibremark}, to conclude that the pull--back 
$p^*$ preserves the curvature types. 
\end{proof}

\medbreak
We need to recall a few facts about {\em basic} connections 
in the foliated Hermitian vector bundle $(E, \npart_E)$ ~\cite{KTone}~. 
Note that all the statements below are of local nature. 

\begin{lem}\label{basic.conn}
For any adapted connection $\nabla_E$~, the following 
conditions are equivalent~$:$ 

\begin{itemize}
\item[$(1)$]
The contraction $i_X \nabla^2_E = 0~,$ for all vector fields 
$X$ in $T(\pi)~;$

\item[$(2)$]
The mixed component $(\nabla_E^2)^{1,1}$ of $\nabla_E$ vanishes$~;$ 

\item[$(3)$]
The curvature $\nabla_E^2$ coincides with the basic component 
$(\nabla_E^2)^{2,0}~,$ that is $\nabla_E^2 = (\nabla_E^2)^{2,0}~;$ 

\item[$(4)$]
For any $\pi$--projectable transversal vector field $\wti{Y},$ 
the operator $\nabla_{\wti{Y}}$ preserves the sheaf $\ker \npart_E$ and 
de\-pends only on $Y = \pi_* \wti{Y}.$ 

\medbreak\noindent
The following condition is a consequence of the above properties$~:$ 

\item[$(5)$]
For $\pi$--projectable transversal vector fields $\wti{Y}, \wti{Y}',$ 
the curvature $\nabla_E^2 (\wti{Y}, \wti{Y}')$ preserves $\ker \npart_E$ 
and de\-pends only on $Y = \pi_* \wti{Y}, Y' = \pi_* \wti{Y}'.$ 
\end{itemize}
\end{lem}

\begin{proof}
Since $(\nabla_E^2)^{0,2} = (\npart_E)^2 = 0~,$ 
the equivalence of $(1), (2)$ and $(3)$ is immediate, 
so we ela\-bo\-rate only on conditions $(4)$ and $(5)$~. 
The mixed component $(\nabla_E^2)^{1,1}$
is characterized by the formula 
\begin{equation}\label{mixed.curv}
(\nabla_E^2)^{1,1} (X, Y) (s) = 
\nabla_E^2 (X, \wti{Y}) (s) = \npart_X (\nabla_{\wti{Y}} s) - 
\nabla_{\wti{Y}} (\npart_X s) - \npart_{[X, \wti{Y}]} s~, 
\end{equation}
for vector fields $X$ in $T(\pi)$ and $\pi$--projectable 
transversal vector fields $\wti{Y}$~. 
Thus for $s \in \ker \npart$~, we have from $(2)~$ 
\begin{equation}\label{parallel1}
(\nabla_E^2)^{1,1} (X, Y) (s) = \npart_X (\nabla_{\wti{Y}} s) = 
\npart_X (\nabla_Y s) \equiv 0~, 
\end{equation}
since $\pi_* [X, \wti{Y}] = [\pi_* X, \pi_* \wti{Y}] = [0, Y] = 0$~. 
The implication $(4) \Ra (2)$ follows from \eqref{parallel1}, 
using Lemma \ref{parallel.sec}. 

\noindent
Likewise, the vector fields $\wti{Y}, \wti{Y}'$ satisfy  
$\pi_* [\wti{Y}, \wti{Y}'] = [\pi_* \wti{Y}, \pi_* \wti{Y}'] = [Y, Y']$ 
and we have from $(4)$ for $s \in \ker \npart$~:
\begin{equation}\label{basic.curv}
\begin{aligned}
(\nabla_E^2)^{2,0} (Y, Y') (s) = \nabla_E^2 (\wti{Y}, \wti{Y}') (s) &= 
\nabla_{\wti{Y}} (\nabla_{\wti{Y}'} s) - 
\nabla_{\wti{Y}'} (\nabla_{\wti{Y}} s) - \nabla_{[{\wti{Y}}, \wti{Y}']} s \\
&= \nabla_Y (\nabla_{Y'} s) - 
\nabla_{Y'} (\nabla_Y s) - \nabla_{[Y, Y']} s~. 
\end{aligned}
\end{equation}
Thus $(4)$ implies $(5)$~. 
\end{proof}

We say that $\nabla_E$ is a {\em basic} connection, if any 
of the equivalent conditions in Lemma \ref{basic.conn} holds. 
In general, a foliated vector bundle $(E, \npart_E)$ does not 
admit basic connections. In our context, the following example 
describes essentially the class of foliated bundles which do 
admit basic connections. 

\begin{ex}\label{loc.triv1}
Locally trivial families of flat bundles~:~
As in Example \ref{holonomy1}, we view a foliated vector bundle 
$(E, \npart_E)$ as a family of flat bundles on the fibers $T_b~,$ 
parametrized by $b \in B~.$ 
We say that this family is {\em locally trivial} if there exists 
a flat bundle $(E_0, \nabla_0)$ on the torus $T~,$ 
determined by a holonomy homomorphism 
$\xi = (\xi_1, \ldots, \xi_n) \in \Hom_{\bZ}(\Lambda, \U(n))~,$ 
and a (good) open cover $\cU$ of $B$ such that for every $U \in \cU~,$ 
there are isomorphisms of foliated vector bundles as indicated in 
the following diagram, similar to \eqref{holonomy2}. 
\begin{equation}\label{loc.triv2}
\begin{CD}
U \times (\bR^d \times_{\Lambda} \bC^n) @>{\cong}>> U \times E_0 
@>{\cong}>>  E \mid \pi^{-1} (U)  \\
@VV{\Id \times \tau'}V  @VV{\Id \times \tau_0}V  
@VV{\tau\vert\pi^{-1}U}V  \\
U \times T    @>{=}>>  U \times T    @>{\cong}>>   \pi^{-1} (U)~. 
\end{CD}
\end{equation}
On overlaps $U_{ik}$ in $\cU~,$ the coordinate changes  
on the LHS are given by \eqref{holonomy4}, except that now 
the holonomy homomorphisms $\xi^i$ are independent of $b \in U_i~.$ 

In the case of locally trivial families, much more can be said 
about the spectral covering $\Sigma$ of $(E, \npart_E)~.$ 
In short, we claim that $\Sigma$ is a finite union of {\em leaves}, 
that is maximal integral manifolds, of the transverse foliation 
$\hcF$ on $\hat\pi:\what{M} \to B$ determined by \eqref{gen.flat2}, 
the leaves of $\hcF$ being holonomy coverings over $B~.$ 
Locally over $U, ~\Sigma_U = \hpis^{-1} (U)$ 
is given by a finite number of constant sections of the 
corresponding trivialization $U \times \what{T} \to U$ of 
$\what{M} \to B$ and we have $\Sigma_{\reg} = \Sigma~,$ 
that is all points of $\Sigma$ are regular and the branch 
locus is empty. 
Therefore $\Sigma \hra \what{M}$ is a smooth submanifold and the 
rank of the multiplicity sheaf $\cL$ is locally constant on 
$\Sigma~,$ hence constant on the connected components of $\Sigma~.$ 
The above local properties of $\Sigma_U$ imply that the 
connected components of the spectral covering $\Sigma$ are 
integral manifolds of the transverse foliation $\hcF~.$ 
The structure of $\Sigma_U$ shows also that 
$\hpis : \Sigma \to B$ satisfies the unique 
path--lifting property. Therefore any path in the leaf 
$\hcF_{\xi}$ through $\xi \in \Sigma~,$ starting at $\xi$ 
must already be in the connected component of $\Sigma$ 
containing $\xi~.$ Thus the connected components of $\Sigma$ 
are maximal integral manifolds of $\hcF~.$ 
These leaves are closed in $\what{M}~,$ since they intersect the 
complete transversals $\what{T}_b$ in $\leq n$ points. 
More precisely, the sets $\Sigma_b \subset \what{T}_b$ 
are invariant under the action of $\pi_1 (B, b)$ 
on $\what{T}_b$ determined by 
$\rho^* : \pi_1 (B) \to \Aut (\what{T}) \cong \Aut (\Lambda)~,$ 
with the orbits and their multiplicities corresponding to the 
component leaves of $\Sigma$ and the rank of $\cL$ on these 
components respectively. This allows us to decompose 
$(\cL, \Sigma)~,$ respectively $(E, \npart_E)$ according to  
its leaf components. 

We abbreviate the above properties of $\Sigma$ by saying that 
$\Sigma \hra \what{M}$ is {\em locally constant}. 
From the local formula \eqref{p-curv} for the curvature $\bF$ 
of the Poincar\'e bundle $\rpb$~, we see that 
$\bF \mid Z_{\Sigma_U} = 0$~. 
Obviously, $\nabla_0$ extends to a basic, in fact a flat connection 
$\nabla_U$ on $U \times E_0$ and these connections can be patched 
together to a basic connection $\nabla_E$ on $E$ via a partition of 
unity on $B$ subordinate to the cover $\cU$. 
A special case of locally trivial families of flat bundles on the 
fibers is of course given by flat bundles $(E, \nabla_E)$ on the 
total space $M$~(compare Section \ref{flat}), 
in which case the bundle is determined by a global 
holonomy homomorphism $\tilde{\rho} : \pi_1 (M) \to \U (n)~,$ 
so that the representation on $\Lambda$ is determined by restriction, 
that is by the diagram 
\begin{equation}\label{holonomy5}
\begin{CD}
\Lambda = \pi_1 (T)  @>{\xi}>>   \U (1)^n  \\
@VV{\iota_*}V    @VV{i}V   \\
\pi_1 (M) @>{\tilde{\rho}}>>  \U (n)~. 
\end{CD}
\end{equation}
\end{ex} 

\medbreak
In order to understand the interplay between the obstruction for 
the existence of a basic connection and the behaviour of the curvature 
term $\bF \mid Z_\Sigma$~, we next look at the case of foliated complex 
line bundles. 

\begin{ex}\label{line1}
Suppose that $(E, \npart_E)$ is a foliated complex line bundle 
on $M$~. In this case, the spectral covering $\Sigma \hra \what{M}$ 
is a section $\sigma$ of $\what{M} \to B$~, the multiplicity 
sheaf $\cL$ on $\Sigma$ is a sheaf of rank $1$ and 
$\ps : Z_\Sigma \to M$ is a diffeomorphism. 
Thus by Theorem \ref{mainresult1} we have 
\begin{equation}\label{line2}
\ps^* E \cong \hps^* \cL \ot \rpb_\Sigma~.
\end{equation} 
It follows that the connection 
$\nabla_\rpb$ and any connection $\nabla_{\cL}$ on $\cL$ induce 
an adapted connection $\nabla_E$ on $E$~, 
such that $\ps^* ((\nabla_E^2)^{1,1}) = \bF \mid Z_\Sigma$ 
(compare Section \ref{inverse2})~. 
This connection $\nabla_E$ is basic, if (and only if) the spectral 
section $\sigma$ is locally constant, that is $(E, \npart_E)$ is a 
locally trivial family of flat line bundles. This follows from 
Example \ref{loc.triv1}. 

For foliated line bundles, the functor $\fm : E \bra \cL$ has the 
following multiplicative property. 
Given $E = E_1 \ot E_2~,$ the spectral sections are 
related by $\sigma = \sigma_1 + \sigma_2$ and we denote by 
$p_i : \Sigma \to \Sigma_i$ the canonical projection. Then a 
direct calculation from 
$\ps^* E \cong \hps^* \cL \ot \rpb_\Sigma~, ~
p_{\Sigma_i}^* E_i \cong\hat{p}_{\Sigma_i}^* \cL_i \ot\rpb_{\Sigma_i}$ 
yields the product formula on $\Sigma~,$ respectively $Z_\Sigma$ 
\begin{equation}\label{line3}
\begin{aligned}
\cL \cong p_1^* ~\cL_1 &\ot p_2^* ~\cL_2~, \\
\rpb_{\Sigma} \cong (\id \times p_1)^*~ \rpb_{\Sigma_1} &\ot 
(\id \times p_2)^* ~\rpb_{\Sigma_2}~. 
\end{aligned}
\end{equation}
\end{ex}

\medbreak
These examples motivate the following definition. 

\begin{dfn}\label{pbasic1}
The adapted connection $\nabla_E$ is {\em Poincar\'e basic} 
if the $\hat{p}_\Sigma$--adapted connection $\wti{\nabla}_E$
on $p_\Sigma^*E \otimes \rpbv_\Sigma$ on $Z_\Sigma$ is basic, 
that is from \eqref{curv2} in Lemma \ref{anticomlemma}
\begin{equation}\label{pbasic2}
\tilde{j}^* (\wti{\nabla}_E^2)^{1,1} = \ps^* (\nabla_E^2)^{1,1} 
- \bF \mid Z_\Sigma = 0~.
\end{equation}
Here we view the scalar form $\bF$ as a form with values in the 
center of $\ps^* \End_s (E)~,$ 
using the canonical isomorphism of foliated bundles 
$\End_s (p^* E \ot \rpbv) \cong p^* \End_s (E)~.$ 

\noindent
For our purposes, it is actually sufficient that the curvature 
$\tilde{j}^* (\wti{\nabla}_E^2)^{1,1}$ vanishes on 
the subsheaf $\cK= \tilde{j}^* (\ker\nabla_E^r)$ defined earlier, 
that is we have $\tilde{j}^* (\wti{\nabla}_E^2)^{1,1} \mid \cK = 0$ 
or equivalently 
\begin{equation}\label{pbasic3}
\ps^* (\nabla_E^2)^{1,1} \mid \cK = 
\bF \mid Z_\Sigma~.
\end{equation}
This corresponds to the equivalent condition $(4)$ in Lemma 
\ref{basic.conn} applied to the connection $\wti{\nabla}_E~.$ 
\end{dfn}

\medbreak
We now proceed to construct a connection
$\nabla_{\cL} : \cL \lra \cL \ot \Omega^1_\Sigma$~, 
given a Poincar\'e basic connection $\nabla_E$ on $E \lra M$~.
From \eqref{curv2} we see that  
\begin{equation}\label{xisigma}
\Xi \mid {Z_\Sigma} = 0 : p_{\Sigma}^*E\otimes\rpbv_{\Sigma}\lra
p_{\Sigma}^*E\otimes\rpbv_{\Sigma}\otimes \Omega^{1,1}_ {Z_\Sigma}~.
\end{equation}
Therefore the diagram below with exact rows is commutative (up to sign) 
\begin{equation}\label{diagram1}
\xymatrix{ 0 \ar[r] & \cK \ar[r] \ar[d] &
\ps^*E \ot \rpbv_{\Sigma} \ar[r]^{\tilde{j}^*\nabla_E^r}
\ar[d]^{\tilde{j}^*\nabla_E^t} &
\ps^*E\ot\rpbv_{\Sigma} \ot \Omega^1_{Z_\Sigma/\Sigma}
\ar[r]^{\tilde{j}^*\nabla_E^r}
\ar[d]^{\tilde{j}^*\nabla_E^t} & 
\ps^*E \ot \rpbv_{\Sigma} \ot \Omega_{Z_\Sigma}^{0,2}
\ar[d]^{\tilde{j}^*\nabla_E^t}  \\ 
0 \ar[r] & \cK \ot \hps^*\Omega^1_\Sigma 
\ar[r] \ar[d] &
\ps^*E\otimes\rpbv_{\Sigma} \ot \hps^*\Omega^1_\Sigma 
\ar[r]^{\tilde{j}^* \nabla_E^r} 
\ar[d]^{\tilde{j}^*\nabla_E^t}  &
\ps^*E \ot \rpbv_{\Sigma} \ot \Omega_{Z_\Sigma}^{1,1} 
\ar[r]^{\tilde{j}^* \nabla_E^r} 
\ar[d]^{\tilde{j}^*\nabla_E^t}  &
\ps^*E \ot \rpbv_{\Sigma} \ot \Omega_{Z_\Sigma}^{1,2} 
\ar[d]^{\tilde{j}^*\nabla_E^t} \\ 
0 \ar[r] & \cK \ot \hat{p}_\Sigma^*\Omega^2_{\Sigma} \ar[r] & 
\ps^* E \ot \rpbv_{\Sigma}\otimes
\hps^*\Omega^2_\Sigma \ar[r]^{\tilde{j}^* \nabla_E^r}  &
\ps^*E \ot \rpbv_{\Sigma} \ot \Omega^{2,1}_{Z_\Sigma} 
\ar[r]^{\tilde{j}^* \nabla_E^r} &
\ps^*E \ot \rpbv_{\Sigma} \ot \Omega^{2,2}_{Z_\Sigma}  }
\end{equation}
and the restriction of $\tilde{j}^*\nabla_E^t$ to the subsheaf $\cK$
induces a connection 
\begin{equation}\label{nablaker}  
\nabla_E^{\ker}: \cK \lra \cK\otimes\hps^*\Omega^1_\Sigma~. 
\end{equation}
Alternatively, we may use the 
condition \eqref{pbasic3} to arrive at the same conclusion. 
Recalling that $\cL = \cR^d \hpsstar (\cK)$ 
and using the projection formula, this leads to a connection 
\begin{equation} \label{transf.conn}
\nabla_{\cL} = \cR^d \hpsstar (\nabla_E^{\ker})~:~
\cL \lra \cL \ot \Omega^1_\Sigma~.
\end{equation}
For later reference, we compute the curvature of the transformed
connection $\nabla_{\cL}$~.

\begin{lem}\label{transf.curv1}
The curvature of the connection $\nabla_{\cL}$ is given by
\begin{equation}\label{transf.curv2}
\nabla_{\cL}^2 = \cR^d  \hpsstar(\nabla_E^{\ker})^2 =
\cR^d  \hpsstar (\ps^* (\nabla_E^2)^{2,0} \mid \cK)~. 
\end{equation}
\end{lem}

\begin{proof}
Since $\tilde{j}^* (\wti{\nabla}_E^2)^{0,2} = 
\tilde{j}^* (\wti{\nabla}_E^2)^{1,1} = 0$ by assumption, 
the curvature term $\tilde{j}^* \wti{\nabla}_E^2 = 
\tilde{j}^* (\wti{\nabla}_E^2)^{2,0} = \ps^* (\nabla_E^2)^{2,0}$ 
leaves the sheaf $\cK \to Z_\Sigma$ invariant by Lemma 
\ref{basic.conn}, $(5)$ and we have from \eqref{curv1} 
$$
\begin{aligned}
\nabla_{\cL}^2 &=  \cR^d \hpsstar (\nabla_E^{\ker})^2 \\ &=
\cR^d \hpsstar (\tilde{j}^*\nabla_E^t \mid {\cK})^2   \\ &=
\cR^d \hpsstar (\tilde{j}^*(\nabla_E^t)^2 \mid {\cK}) \\ &=
\cR^d \hpsstar (\tilde{j}^* (p^* (\nabla_E^2)^{2,0}) \mid {\cK}) \\ &=
\cR^d \hpsstar (\ps^* (\nabla_E^2)^{2,0} \mid {\cK})~.
\end{aligned}
$$
In this calculation we also used diagram \eqref{diagram1}. 
\end{proof}

\medbreak
Recall that the pair $(E,\nabla_E)$ is said to be {\it reducible}
if there are bundles with connections $(E_1,\nabla_{E_1})$ and
$(E_2,\nabla_{E_2})$ such that $E = E_1 \oplus E_2$ and
$\nabla_E = \nabla_{E_1} \oplus \nabla_{E_2}$~. The pair
$(E,\nabla_E)$ is said to be {\it irreducible} if it is not
reducible.

\begin{lem} \label{sum1}
If $(E, \nabla_E) = (E_1, \nabla_{E_1}) \oplus (E_2,\nabla_{E_2})$, 
then
\begin{equation}
( \cL,\nabla_{\cL}, \Sigma ) = 
( \cL_1, \nabla_{\cL_1}, \Sigma_1 ) \oplus 
( \cL_2, \nabla_{\cL_2}, \Sigma_2)~, 
\end{equation}
where $\Sigma = \Sigma_1 \cup \Sigma_2~.$ 
\end{lem}

\begin{proof}
The statement is clear from the definitions. Indeed, we have
$\ker\nabla_E = \ker\nabla_{E_1}\oplus\ker\nabla_{E_2}$, so that
$\what E = \what E_1 \oplus \what E_2$ by (\ref{transf.bdl1}) and
therefore $\cL = \cL_1 \oplus \cL_2$~. Since $\nabla_E^{\ker}$
also splits as a direct sum, it follows from \eqref{transf.conn}
that $\nabla_{\cL} = \nabla_{\cL_1} \oplus \nabla_{\cL_2}$~.
As for supports, we note that $\Sigma = \Sigma (E, \npart_E) = 
\Sigma_1 \cup \Sigma_2~,$ where $\Sigma_i = \Sigma (E_i, \npart_{E_i})~.$ 
\end{proof}

\begin{dfn}
The triple $(\cL, \nabla_{\cL}, \Sigma)$ is called the 
{\it Fourier--Mukai transform} of $(E,\nabla_E)$, and 
$\Sigma= \Sigma (E, \npart_E) = \supp \what{E}$ is called the 
{\it spectral covering} associated with the underlying 
foliated structure $(E,\npart_E)$~.
\end{dfn}

\subsection{The inverse transform for connections}
\label{inverse2}

Given $(\cS, \Sigma)$ as in Section \ref{inverse1}, 
let $\nabla_{\cS}$ be a connection on $\cS~.$ 
In order to obtain a connection on $\fmh (\cS) = \check{\cS}~,$ 
recall from \eqref{checks} that $\check{\cS}$ is defined by 
$\check{\cS} = \psstar (\hps^* \cS \ot \rpb_\Sigma)~.$ 
Consider the connection 
\begin{equation} \label{tilde.conn2}
\wti{\nabla}_{\cS} = \hps^*\nabla_{\cS}
\ot \id_{\rpb_\Sigma} + \id_{\cS} \ot \nabla_{\rpb_\Sigma} :
\hps^*\cS \ot \rpb_\Sigma \lra \hps^* \cS \ot 
\rpb_\Sigma \ot \Omega^1_{Z_\Sigma}~.
\end{equation}
Since $\Omega^1_{Z_\Sigma}=\ps^*\Omega^1_M$, it follows from the 
projection formula that
\begin{equation}
\psstar (\hps^*\cS \ot \rpb_\Sigma \ot \Omega^1_{Z_\Sigma}) = 
\check{\cS} \ot \Omega^1_M~.
\end{equation}
Thus we define the connection $\nabla_{\check{\cS}}$ on $\check{\cS}~,$ 
adapted to the flat partial connection $\npart_{\check{\cS}}$ in 
Section \ref{inverse1}, by 
\begin{equation} \label{transf.conn2}
\nabla_{\check{\cS}} = \psstar \wti{\nabla}_{\cS}~: ~\check{\cS}
\lra \check{\cS}\otimes\Omega^1_M~.
\end{equation}
The curvature $\nabla_{\check{\cS}}^2$ of $\nabla_{\check{\cS}}$ 
is computed next, from the formula    
\begin{equation}\label{inv.curv}
\wti{\nabla}_{\cS}^2 = \hps^*\nabla_{\cS}^2 
\otimes\id_{\rpb_\Sigma} + \id_{\cS} \otimes \bF \mid Z_\Sigma~. 
\end{equation}

\begin{lem} \label{inv.curv1}
The curvature $\nabla_{\check{\cS}}^2$ of the connection 
$\nabla_{\check{\cS}}$ is determined by 
\begin{equation}\label{inv.curv2}
(\nabla_{\check{\cS}}^2)^{0,2} = 0  \qquad,\qquad 
(\nabla_{\check{\cS}}^2)^{2,0} = \psstar (\hps^* \nabla_{\cS}^2)~, 
\end{equation}
and 
\begin{equation}\label{inv.curv3}
(\nabla_{\check{\cS}}^2)^{1,1} = \psstar (\nabla_{\rpb_\Sigma}^2) = 
\psstar (\bF \mid Z_\Sigma)~.
\end{equation}
This implies that the connection $\nabla_{\check{\cS}}$ is 
Poincar\'e basic, since the pull--back $\hps^*\nabla_{\cS}$ 
is $\hps$--basic. 
\end{lem}

\medbreak
The triple $(\cS,\nabla_{\cS}, \Sigma)$ is reducible if there are triples 
$(\cS_1, \nabla_{\cS_1}, \Sigma_1)$ and $(\cS_2, \nabla_{\cS_2}, \Sigma_2)$, 
such that we have $\Sigma = \Sigma_1 \cup \Sigma_2$ and 
$(\cS, \nabla_{\cS}, \Sigma) = (\cS_1, \nabla_{\cS_1}, \Sigma_1) 
\oplus (\cS_2, \nabla_{\cS_2}, \Sigma_2)$~.
The triple $(\cS, \nabla_{\cS}, \Sigma)$ is said to be irreducible 
if it is not reducible. 
Moreover, a $(m-d)$--dimensional smooth submanifold 
$\Sigma \hra \what M~,$ which is an $n$--fold covering of $B$ 
transversal to all fibers, is said to be {\it proper} if the trivial 
local system $(\underline{\bC},\trivcon, \Sigma)~,$ consisting of the 
trivial line bundle on $\Sigma$ with the trivial flat connection, 
is irreducible. 
Clearly, if $\Sigma$ is proper, then it is connected. 

\begin{lem} \label{sum2}
If
$( \cS, \nabla_{\cS}, \Sigma) = 
( \cS_1, \nabla_{\cS_1}, \Sigma_1 ) \oplus 
( \cS_2, \nabla_{\cS_2}, \Sigma_2 )~,$ 
then
\begin{equation}
( \check{\cS}, \nabla_{\check{\cS}}) = 
(\check{\cS}_1, \nabla_{\check{\cS}_1}) \oplus 
(\check{\cS}_2, \nabla_{\check{\cS}_2})~.
\end{equation}
Moreover, if $(\cS,\nabla_{\cS}, \Sigma)$ is irreducible 
with smooth support $\Sigma~,$ 
then $\Sigma = \supp \cS$ is proper.
\end{lem}

\begin{proof}
The first statement follows easily from the definitions of the
previous paragraph.
Now if $\Sigma = \supp \cS$ is not proper, then 
$(\underline{\bC},\trivcon, \Sigma)$ splits as the sum 
$(\cS_1, \nabla_{\cS_1}, \Sigma_1) \oplus 
(\cS_2,\nabla_{\cS_2}, \Sigma_2)~,$ where $\Sigma = 
\Sigma_1 \cup \Sigma_2~, ~\Sigma_i = \supp \cS_i~.$ 
Thus
$$ 
\begin{aligned}
(\cS, \nabla_{\cS}, \Sigma) &= (\cS, \nabla_{\cS}, \Sigma) 
\ot (\underline{\bC}, \trivcon, \Sigma)  \\ 
&= (\cS, \nabla_{\cS}, \Sigma) \ot (\cS_1,\nabla_{\cS_1}, \Sigma_1) 
\oplus (\cS, \nabla_{\cS}, \Sigma) \ot (\cS_2, \nabla_{\cS_2}, \Sigma_2) \\ 
&= (\cS \ot \cS_1, \nabla_{\cS \ot \cS_1}, \Sigma_1) \oplus 
(\cS \ot \cS_2, \nabla_{\cS \ot \cS_2}, \Sigma_2)~,
\end{aligned}
$$
and $(\cS, \nabla_{\cS}, \Sigma)$ is reducible. 
\end{proof}

\subsection{The main theorem for bundles with Poincar\'e 
basic connections} \label{main2}

\begin{dfn}
$\Vect{n}{M}$ is the category of foliated Hermitian vector bundles 
on $M$ endowed with a Poincar\'e basic unitary connection. 
Objects in $\vectn$ are pairs $(E,\nabla_E)$ consisting of a 
foliated Hermitian vector bundle $E$ of rank
$n$ and a Poincar\'e basic unitary connection $\nabla_E$~. 
Morphisms are bundle maps compatible with the connections.
\end{dfn}

\begin{dfn}
$\Spec{n}$ is the category of {\it spectral data} on $\what M$.
Objects in $\Spec{n}$ are triples $(\cS,\nabla_{\cS}, \Sigma)~,$ 
such that the pair $(\cS, \Sigma)$ is an object in $\relskyn$ and 
$\nabla_{\cS}$ is a connection on $\cS\vert_{\Sigma}$~.
Morphisms are sheaf maps of $\cO_{\what{M}}$--modules 
compatible with the connections.
\end{dfn}

\begin{thm} \label{mainresult2}
The Fourier--Mukai trans\-form $\fm$ defines an additive 
natural equi\-valence of cate\-gories
\begin{equation}
\fm~:~ \Vect{}{M} \osetl{\cong}  \Spec{}~.
\end{equation}
\end{thm}

\begin{proof}
In view of the natural isomorphisms 
$\phi_{\cS} : \cS \oset{\cong} \fm \circ \fmh (\cS)~, ~
\psi_E : \fmh \circ \fm (E) \oset{\cong} E$ 
in the proof of Theorem \ref{mainresult1}, it suffices to show that 
we have gauge equivalences 
$\phi_{\cS}^* \nabla_{\cL_{\check{\cS}}} = \nabla_{\cS}$ and 
$\psi_E^* \nabla_E = \nabla_{\check{\cL}}~.$
We comment only on the proof for the second gauge equivalence. 
In fact, our constructions of $\fm$ in Section \ref{basic1} and 
$\fmh$ in Section \ref{inverse2} show that $\nabla_E$ corresponds 
to $\psstar (\hps^* \wti{\nabla}_{\cL})~.$ 
\end{proof}

\medbreak
Let $\Vect{n}{B}$ be the category of complex vector bundles $V$ 
of rank $n$ over $B$ with unitary connection $\nabla_V$ and fix 
a pair $(V,\nabla_V) \in \Vect{n}{B}~.$ 
Then $\pi^*(V, \nabla_V) = (\pi^* V, \pi^* \nabla_V)$ is an 
object in $\Vect{n}{M}$, while $\hat{\pi}_0^* (V, \nabla_V)$ 
is an object in $\Spec{n}$ supported on the $0$--section 
$\Sigma_0 = \sigma_0 (B) \subset \what{M}~.$ 
The construction of $\fm$ is again compatible with these 
pull--backs, that is we have the commutative diagram similar to 
\eqref{diagram3} 
\begin{equation}\label{diagram4}
\begin{CD}
\Vect{n}{M} @>{\fm}>> \Spec{n}      \\
@A{\pi^*}AA   @A{\hat\pi_0^*}AA   \\
\Vect{n}{B}  @>{=}>> \Vect{n}{B} . 
\end{CD}
\end{equation}
Moreover, Corollary \ref{module1} remains valid for 
$\fm$ on $\Vect{}{M}~.$

\begin{cor}\label{module3}
For $(E,\nabla_E) \in \Vect{}{M}$ and $(V,\nabla_V) \in \Vect{}{B}~,$ 
the Fourier--Mukai trans\-form $\fm$ satisfies 
\begin{equation}\label{module4}
\fm (\pi^* (V,\nabla_V) \ot (E,\nabla_E)) \cong 
\hpis^* (V,\nabla_V) \ot \fm (E,\nabla_E)~,
\end{equation}
where $\Sigma$ is the support of $\fm (E,\nabla_E)~.$ 
\end{cor}

\begin{cor}\label{pullback1}
$(E,\nabla_E) \in \Vect{n}{M}$ is of the form 
$(E,\nabla_E) = \pi^* (V,\nabla_V)~,$ for 
$(V,\nabla_V) \in \Vect{n}{B}~,$ if and only if the support 
of the Fourier--Mukai trans\-form $\fm (E,\nabla_E)$ 
is the $0$--section 
$\Sigma_0 = \sigma_0 (B)$ of $\hat{p} : \what{M} \to B~.$ 
\end{cor}

\medbreak
As a consequence of Lemmas \ref{sum1}, \ref{sum2} and 
Theorem \ref{mainresult2}, we have~:

\begin{cor}\label{irred1}
The pair $(E,\nabla_E)$ is irreducible, if and only if its
transform $\fm(E,\nabla_E)$ is irreducible.
\end{cor}

\medbreak
For complex vector bundles $(E, \nabla_E)$ with unitary 
connection $\nabla_E~,$ there is a well-known reduction theorem 
\cite{KN}, Ch. II,~Thm.~7.1~, 
based on the decomposition of the holonomy group in $\U (n)$ 
into irreducible components. Our construction shows that this 
defines a decomposition of $(E, \nabla_E) \in \Vect{}{M}$ 
into irreducible components. From Theorem \ref{mainresult2} 
and Corollary \ref{irred1}, 
we obtain a similar decomposition of 
$\fm (E, \nabla_E)$ in $\Spec{}~.$ 
In the smooth case, the irreducible pairs 
$(E, \nabla_E) \in \Vect{n}{M}$ are characterized as follows. 

\begin{prop}\label{irred2}
Suppose that the spectral covering $\Sigma$ of 
$(E, \nabla_E) \in \Vect{n}{M}$ 
is smooth. Then the pair $(E, \nabla_E)$ is irreducible, if and only if 
its transform $(\cL,\nabla_{\cL}, \Sigma) = \fm (E, \nabla_E)$ 
satisfies the following conditions$~:$ 
$\Sigma$ is connected, $\abs{\Sigma_b} \equiv \ell$ on $B$ 
for some $\ell \vert n~,$ and $(\cL, \nabla_{\cL})$ is a 
vector bundle of rank $k = \frac{n}{\ell}$ with irreducible 
holonomy. 
The covering  $\hat{\pi}_{\Sigma}:\Sigma \to B$ is non--degenerate exactly 
for $\ell = n~, ~k = 1~.$ 
In addition, any smooth, connected spectral manifold 
$\Sigma \hra \what{M}$ is proper. 
\end{prop}

\begin{proof}
This is a consequence of Corollary \ref{irred1} and the remarks 
following Proposition \ref{two}.
\end{proof}

\medbreak
From Example \ref{loc.triv1}, we have the following 
characterization of locally trivial families $(E, \npart_E)$ 
of flat bundles along the fibers. 

\begin{cor}\label{loc.triv3}
The Fourier--Mukai trans\-form $\fm$ defines an equivalence between  
pairs $(E,\nabla_E) \in \Vect{n}{M}~,$ such that $(E, \npart_E)$ is a 
locally trivial family of flat bundles along the fibers and 
$\nabla_E$ is basic$;$ and spectral data 
$(\cS,\nabla_{\cS}, \Sigma) \in \Spec{n}~,$ such that the 
spectral covering $\Sigma \subset \what{M}$ is locally constant 
and $\cS$ has locally constant rank on $\Sigma~.$ 
\end{cor}

Combining this with Proposition \ref{irred2}, we obtain in addition 

\begin{cor}\label{irred3}
$(E,\nabla_E) \in \Vect{n}{M}$ as in Corollary $\ref{loc.triv3}$ 
is irreducible, if and only if the spectral covering 
$\Sigma$ is a leaf of the transversal foliation 
$\hcF$ in \eqref{gen.flat2}, $\abs{\Sigma_b} \equiv \ell$ 
on $B$ for some $\ell \vert n~, ~\pi_1 (B, b)$ acts transitively on 
$\Sigma_b = \Sigma \cap \what{T}_b$ and $(\cL, \nabla_{\cL})$ is a 
vector bundle of rank $k = \frac{n}{\ell}$ with irreducible holonomy. 
\end{cor}

\medbreak
From Example \ref{line1}, in particular formula \eqref{line3}, 
we have the following characterization of foliated line bundles on $M~.$ 

\begin{cor}\label{line5}
The Fourier--Mukai trans\-form $\fm$ defines a multi\-pli\-ca\-tive 
equi\-va\-lence bet\-ween pairs $(E,\nabla_E) \in \Vect{1}{M}$ 
and spectral data $(\cS,\nabla_{\cS}, \Sigma) \in \Spec{1}~,$ 
where $\cS$ is a complex line bundle with connection $\nabla_{\cS}$ 
on the spectral section $\Sigma = \sigma (B)~.$ 
\end{cor}


\section{Applications and examples} \label{examples}

Let us now apply our theorems to a few interesting examples.
We are particularly interested in seeing how differential conditions
on the connection $\nabla_E$ are transformed.

\subsection{Local systems and the representation variety} \label{flat}

As a first example, we now look at the action of the Fourier--Mukai
transform in the subcategory $\locmn$ of unitary local systems 
of rank $n$ on $M~.$ 
So let $E \lra M$ be a complex Hermitian vector bundle of rank $n$ 
and take $\nabla_E$ to be a flat unitary connection on $E$~.

\begin{lem} \label{flat1}
If $\nabla_E$ is flat, then the spectral covering 
$\Sigma \hra \what{M}$ of $(E, \npart_E)$  is locally 
constant, the rank of $\cL$ is locally constant on $\Sigma$ 
and the transform $\nabla_{\cL}$ is also flat.
\end{lem}

\begin{proof}
This follows from \eqref{transf.curv2}, observing that 
a flat connection is basic with locally constant spectral covering 
$\Sigma \hra \what{M}$ and hence $\bF \mid Z_{\Sigma} =0~$ 
(compare Example \ref{loc.triv1}).
\end{proof}

Next, we argue that the inverse transform also preserves flatness, 
provided the spectral covering $\Sigma \hra \what{M}$ is 
locally constant.

\begin{lem}\label{flat2}
If the spectral covering $\Sigma \hra \what{M}$ is locally 
constant and $\nabla_{\cS}$ is flat, then its transform 
$\nabla_{\check{\cS}}$ is also flat. 
\end{lem}

\begin{proof}
Firstly, the assumption on the spectral covering implies 
that $\nabla^2_{\rpb} \mid Z_{\Sigma} = \bF \mid Z_{\Sigma} = 0~.$ 
Since $\nabla_{\cS}$ is flat, the Lemma follows from 
\eqref{inv.curv} or Lemma \ref{inv.curv1}. 
\end{proof}

With these facts in mind, we introduce the following definition:

\begin{dfn}
$\Lochmn{n}$ is the full subcategory of $\specn$ consisting of those
objects $(\cS,\nabla_{\cS}, \Sigma)$ such that the spectral covering 
$\Sigma \hra \what{M}$ and the rank of $\cS$ on $\Sigma$ are locally 
constant, and $\nabla_{\cS}$ is flat. 
\end{dfn}

As a consequence of Theorem \ref{mainresult2} and Lemmas
\ref{flat1} and \ref{flat2}, we obtain~:

\begin{thm}\label{mainresult3}
The Fourier--Mukai transform $\fm$ defines a 
natural equivalence of categories 
$$
\fm~:~\locmn \osetl{\cong} \lochmn~.
$$
\end{thm}

\medbreak
For unitary local systems on $M~,$ the decomposition into 
irreducible components in Section \ref{main2} applies 
mutatis mutandis and we may sharpen Corollary \ref{irred3} 
accordingly, using Theorem \ref{mainresult3}. 

\begin{cor}\label{irred4}
$(E,\nabla_E) \in \locmn$ is irreducible, if and only if the 
spectral covering $\Sigma$ is a leaf of the transversal foliation 
$\hcF$ in \eqref{gen.flat2}, $\abs{\Sigma_b} \equiv \ell$ 
on $B$ for some $\ell \vert n~, ~\pi_1 (B, b)$ acts transitively on 
$\Sigma_b$ and $(\cL, \nabla_{\cL})$ is an irreducible flat vector 
bundle of rank $k = \frac{n}{\ell}~,$ 
that is an irreducible $\U (k)$--local system on $\Sigma~.$ 
\end{cor}


\medbreak
Now let $\cR_M (n)$ denote the moduli space of irreducible unitary
local systems of rank $n \geq 1$ on $M$~. Recall that $\cR_M (n)$
coincides with the {\it representation variety} of $\pi_1(M)$,
that is, the set of all irreducible representations $\pi_1(M) \ra \U(n)$
modulo conjugation. 
Let also $\cS(n)$ denote the set of all connected, locally constant 
$(m-d)$--dimensional submanifolds $\Sigma \hra \what M$ 
(modulo isomorphisms), 
such that the trivial local system $(\underline{\bC},\trivcon,\Sigma)$ 
is a relative skyscraper of length $\ell$ for $\ell \vert n~,$ that is 
\begin{equation}\label{generic1}
(\underline{\bC}, \trivcon, \Sigma) \in \Lochmn{\ell}~. 
\end{equation}
Then \eqref{generic1} implies that $\abs{\Sigma_b} \equiv \ell$ on $B~,$ 
since the multiplicity of this system is $1~.$ 
From Example \ref{loc.triv1} and Corollary \ref{irred4}, we see that 
$\Sigma \hra \what{M}$ is a (proper) $\ell$--sheeted leaf of the 
transverse foliation $\hcF$ on $\hat\pi:\what{M} \to B$ determined by 
\eqref{gen.flat2}, with transitive transversal holonomy of order $\ell~.$ 
For given $\ell \vert n~,$ we denote the corresponding 
subset of $\cS (n)$ by $\cS (n)^{\ell}~.$ 
We will see that the {\em generic} elements of $\cS(n)$ are those 
in $\cS(n)^n$ and we proceed with an explicit parametrization 
of these moduli spaces. 

\medbreak
The transverse foliation $\hcF$ on $\hat\pi: \what{M} \to B$ provides 
the link between our geometric setup and the representation theory and 
we refer again to Example \ref{loc.triv1} and Corollary \ref{irred4} 
for the discussion to follow. 
Recall that the leaves of $\hcF$ are the images 
$\hcF_{\xi}$ of the level sets $\wti{B} \times \{\xi\}$ in 
\eqref{gen.flat2} and are therefore covering spaces over $B$ of 
the form $\hcF_{\xi} \cong \wti{B} / \Gamma_{\rho, \xi}~,$ where 
$\Gamma_{\rho, \xi} \subset \pi_{1} (B, b)$ is the isotropy 
group at $\xi \in \what{T}_b$ under the action corresponding 
to $\rho^* : \pi_1 (B) \to \Aut (\what{T})~.$ 
Here we fix a basepoint $b \in B$ once and for all. 
It is now clear that the structure of the spaces $\cS (n)~,$ 
respectively $\cS (n)^{\ell}$ may be described in terms of the 
transversal holonomy groupoid on the complete transversal 
$\what{T}$ and the leaf space of $\hcF~.$ The leaf space of 
$\hcF$ is the quotient $\pi_1 (B) \bsl \what{T}~,$ 
which may behave quite badly. 
But for our purposes, we need only consider the invariant 
subspace $\what{T}_{\fin} \subset \what{T}$ defined by the 
$\xi \in \what{T}$ satisfying 
$[\Gamma_{\rho, \xi}:\pi_1 (B)] < \infty~,$ that is the 
leaves with finite transversal holonomy, on which the 
$\pi_1 (B) $--orbits are finite by definition. 
$\what{T}_{\fin}$ has an invariant relatively 
closed {\em stratification} 
$\what{T}_{n-1} \subset \what{T}_{n} \subset \ldots~,$ 
given by the points $\xi \in \what{T}_{\fin}$ satisfying 
$[\Gamma_{\rho, \xi} : \pi_1 (B)] \leq n~.$ The main stratum 
in $\what{T}_{n}$ is then given by the invariant relatively 
open set $\what{T}^{n} \subset \what{T}_{n}$ of those points 
$\xi$ for which $[\Gamma_{\rho, \xi}:\pi_1 (B)] = n~.$ 
Here we have tacitly used the `semicontinuity' for the isotropy 
groups of a smooth group action, that is 
$[\Gamma_{\rho, \xi}:\pi_1 (B)] \geq n$ is an open condition and hence 
$[\Gamma_{\rho, \xi}:\pi_1 (B)] \leq n$ is a closed condition. 

For the generic $\Sigma \in \cS (n)^{n}~,$ we have 
$\Sigma = \hcF_{\xi} \cong \wti{B} / \Gamma_{\rho, \xi}~,$ 
where the isotropy group $\Gamma_{\rho, \xi} \subset \pi_{1} (B, b)$ 
has index $n$ and $\Sigma_b = \Sigma \cap \what{T}_b$ corresponds 
to an orbit (of order $n$) in the main stratum $\what{T}^n~.$ 
For $\ell < n~, ~\ell \vert n$ and $\Sigma \in \cS (n)^{\ell}~,$ 
we have $\Sigma = \hcF_{\xi} \cong \wti{B} / \Gamma_{\rho, \xi}~,$ 
where the isotropy group $\Gamma_{\rho, \xi} \subset \pi_{1} (B, b)$ 
has index $\ell$ and $\Sigma_b$ corresponds to the orbit (of order $\ell$) 
of a limit element of $\what{T}^n$ in 
$\what{T}^{\ell} \subset \what{T}_{n-1} = 
\what{T}_{n} \bsl \what{T}^{n}~.$ 
 
\begin{thm}\label{mainresult4}
The space $\cS (n)$ of spectral manifolds associated to irreducible 
$\U (n)$--re\-pre\-sen\-ta\-tions in $\cR_M (n)$ is of the form 
\begin{equation}\label{generic2} 
\cS (n) = \cS (n)^n ~\cup ~( ~\bigcup_{\ell < n}^{\ell \vert n} ~
\cS (n)^{\ell} ~)~, 
\end{equation} 
and is parametrized $($up to automorphisms$)$ by 
\begin{equation}\label{generic3}
\pi_1 (B) ~\bsl ~\{ ~\what{T}^n ~\cup ~
( ~\bigcup_{\ell < n}^{\ell \vert n} ~\what{T}^{\ell} ~) ~\} 
\subset \pi_1 (B) \bsl \what{T}_n~, 
\end{equation}
that is the space of leaves of $\hcF$ with finite transversal holonomy 
of order $\ell$ with $\ell \vert n~.$ 
The mapping 
\begin{equation}\label{generic4}
\Psi (n) : \cR_M (n) \lra \cS (n)~,
\end{equation} 
defined by $\Psi (n) (E,\nabla_E) = {\supp}~\fm(E,\nabla_E) = \Sigma~,$ 
has the following properties$~:$ 

\begin{itemize}
\item[$(1)$]
The generic part $\Psi (n)^n = \Psi (n) \mid \cR_M (n)^n : 
\cR_M (n)^n \to \cS (n)^n$ is surjective, where 
$\cR_M (n)^n = \Psi (n)^{-1} \cS (n)^n$ is 
the space of irreducible representations for which 
the induced fiber holonomy representations 
$\{\xi_1, \ldots, \xi_n\}~, ~\xi_j \in \what{T} \cong 
\Hom_{\bZ}(\Lambda, \U(1))$ of $\Lambda$ consist of $n$ 
distinct elements. The fiber $\Psi (n)^{-1} (\Sigma)$ over 
the generic elements $\Sigma \in \cS(n)^{n}$ corresponds 
exactly to the $\U(1)$--local systems on $\Sigma$ under 
the functor $\fm~.$ 

\item[$(2)$]
The fiber $\Psi (n)^{-1} (\Sigma)$ for $\Sigma \in \cS(n)^{\ell}~, ~
\ell < n~,$ corresponds to $($equi\-valence classes of$)$ 
ir\-redu\-cible $\U(k)$--local systems on $\Sigma$ for 
$k = \frac{n}{\ell}$ under the functor $\fm~.$ 

\end{itemize}
\end{thm}

\medbreak
An extreme situation occurs in the context of Corollary \ref{pullback1}. 
In this case, the irreducible flat vector bundle 
$(E,\nabla_E)$ is the pull--back of an irreducible flat vector 
bundle $(V,\nabla_V)$ on the base $B~,$ determined by an irreducible 
$\U (n)$--representation of $\pi_1 (B)~.$ We have $\ell = 1~, ~k = n~,$ 
$\Sigma$ is the $0$--section $\Sigma_0 = \sigma_0 (B)$ of 
$\hat{\pi} : \what{M} \to B$ and the corresponding orbit is 
given by the origin $0 \in \what{T}^1~.$ 

\begin{proof}
The structure of $\cS (n)$ follows from the description 
preceding the theorem. 
From Corollary \ref{irred4} we see that 
$\Psi (n) (E,\nabla_E) = {\supp}~\fm(E,\nabla_E)$ belongs to $\cS(n)~.$ 
To see that $\Psi (n)^n$ is surjective on $\cR_M (n)^n~,$ 
take $\Sigma\in\cS(n)^n$ and let $(\underline{\bC},\trivcon)$ 
be the trivial $\U(1)$--local
system on $\Sigma$~. Since $\Sigma$ is proper, 
$(\underline{\bC}, \trivcon, \Sigma)$ is irreducible and 
$(\underline{\bC},\trivcon,\Sigma) \in \Lochmn{n}$ implies 
that $\Sigma \in \cS(n)^n$ satisfies $\abs{\Sigma_b} \equiv n~.$ 
Then it follows from Corollary \ref{irred1} and Theorem 
\ref{mainresult3} that $\fmh (\underline{\bC}, \trivcon, \Sigma)$ 
is an irreducible $\U(n)$--local system on $M$ and defines a 
point in $\cR_M (n)^n$ mapping to $\Sigma$ under $\Psi (n)^n~.$ 
The statements about the fibers of $\Psi (n)$ in $(1)$ and $(2)$ 
also follow from Corollary \ref{irred4}, that is 
$(E,\nabla_E) \in \Psi (n)^{-1} (\Sigma)~,$ if and only if 
$\fm(E,\nabla_E)$ is an irreducible $\U(k)$--local 
system on $\Sigma \in \cS (n)^{\ell}~,$ for 
$\ell \vert n~, ~k = \frac{n}{\ell}~.$ 
\end{proof}

\medbreak
For the corresponding irreducible $\U (n)$--representations, 
we recall that the induced fiber holonomy 
representations $\{\xi_1, \ldots, \xi_n\}~, ~\xi_j \in \what{T} \cong 
\Hom_{\bZ}(\Lambda, \U(1))$ of $\Lambda$ consist generically of $n$ 
distinct elements, that is they are orbits of order $n$ of the action 
$\rho^* : \pi_1 (B) \to \Aut (\what{T}) \cong \Aut (\Lambda)~.$ 
How does one then describe the irreducible $\U (n)$--representation 
of $\pi_1 (M)$ in $\Psi (n)^{-1}(\Sigma) \subset \cR_M (n)^n$ associated 
to a $\U(1)$--local system $(\cS, \nabla_{\cS}, \Sigma)$ for 
$\Sigma\in \cS (n)^n~,$ or more generally to any irreducible 
$\U(k)$--local system $(\cS, \nabla_{\cS}, \Sigma)$ for 
$\Sigma \in \cS (n)^{\ell}$ via $\fmh~?$ 
We claim that the irreducible $\U (n)$--representations of the 
crossed product \eqref{cross.prod1} for $\pi_1 (M)$ are obtained 
by the {\em induced representation} from irreducible 
$\U (k)$--representations on a subgroup of index $\ell$ with 
$\ell \vert n$ in $\pi_1 (M)$ to the full group. 
In fact for $k = \frac{n}{\ell}~,$ any $\U(k)$--local system $\cS$ 
on $\Sigma$ is determined by a homomorphism 
$\eta : \pi_1 (\Sigma, \xi) \cong \Gamma_{\rho, \xi} \to \U (k)~,$ 
which together with the above datum $\xi \in \what{T}^{\ell}$ 
defines an irreducible unitary representation 
\begin{equation}\label{induced1}
(\exp (\xi), \eta) : \Lambda \times_{\rho_*} \Gamma_{\rho, \xi} 
\to \U (k)~,
\end{equation}
that is a $\U (k)$--representation on a subgroup of index $\ell$ in 
the crossed product $\pi_1 (M) = \Lambda \times_{\rho_*} \pi_1 (B)~.$ 
We need to verify 
$
\exp (\xi) (\rho_* (\gamma) (a)) = 
\eta (\gamma) \exp (\xi) (a) \eta (\gamma)^{-1} = \exp (\xi) (a)~, ~
a \in \Lambda~, ~\gamma \in \Gamma_{\rho, \xi}~, 
$
which is obvious, since $\U (1)$ is identified with the center of $\U (k)$ 
and $\Gamma_{\rho, \xi}$ fixes $\xi$ under $\rho^*~.$ 

Looking at the construction of the inverse Fourier--Mukai transform 
$\fmh~,$ in particular the push--down operation $\psstar$ for the 
$\ell$--fold covering map $\ps : Z_\Sigma \to M~,$ we see that the 
irreducible $\U (n)$--representation of 
$\pi_1 (M)$ in Theorem \ref{mainresult4}, given by 
$\fmh (\cS_{\eta}, \nabla_{{\cS}_{\eta}}, \Sigma)$ of the 
$\U(k)$--local system $(\cS_{\eta}, \nabla_{{\cS}_{\eta}})$ on 
$\Sigma$ defined by $\eta~,$ is in fact the induced representation 
of $(\exp (\xi), \eta)~.$ 

For the generic case $\ell = n~, ~k = 1~,$ that is 
$(\underline{\bC}, \trivcon, \Sigma) \in \Lochmn{n}~,$ 
the index of $\Gamma_{\rho, \xi} \subset \pi_1 (B)$ is $n$ 
and the unitary local systems on $\Sigma$ are $1$--dimensional. 
It then follows also that the irreducible $\U (n)$--representation 
of $\pi_1 (M)$ associated to 
$(E_{\Sigma}, \nabla_{E_\Sigma}) = \fmh (\underline{\bC}, \trivcon, \Sigma)$
corresponds to the induced representation of $(\exp (\xi), \id)~,$ 
that is the trivial representation $\eta = \id$ of $\Gamma_{\rho, \xi}~.$ 

\medbreak
In summary, the following theorem is the algebraic version of 
Theorem \ref{mainresult4}, given purely in terms of representation 
theory. 

\begin{thm}\label{mainresult4a}
For the torus bundle $\pi : M \to B$ in 
\eqref{torus.bundle1}, the representation variety $\cR_M (n)$ 
of the fundamental group $\pi_1 (M)~,$ given by the crossed 
product \eqref{cross.prod1} with respect to the action 
$\rho : \pi_1 (B) \to \Aut (T) \cong \GL(d, \bZ)~,$ 
is parametrized by the following data$~:$

\begin{itemize}
\item[$(1)$]
Elements $[\xi] \in \pi_1 (B) \bsl \what{T}^{\ell}~,$ for 
$\ell \vert n~,$ that is orbits of order $\ell$ in the dual torus 
$\what{T}$ under the induced action 
$\rho^* : \pi_1 (B) \to \Aut (\what{T})~.$ 

\item[$(2)$]
Irreducible unitary representations 
$\eta : \Gamma_{\rho, \xi} \to \U (k)$ of the isotropy group 
$\Gamma_{\rho, \xi} \subset \pi_1 (B)$ of index $\ell$ 
at $\xi \in \what{T}^{\ell}~,$ for $k = \frac{n}{\ell}~.$ 
\end{itemize}
These data determine an irreducible $\U (n)$--representation 
of $\pi_1 (M)$ by the induced representation of 
$(\exp (\xi), \eta) : 
\Lambda \times_{\rho_*} \Gamma_{\rho, \xi} \to \U (k)$ 
from the subgroup of index $\ell$ in the crossed product $\pi_1 (M)~.$ 
This induction corresponds to the functor $\fmh~.$ 
The generic case occurs for 
$\xi \in \pi_1 (B) \bsl \what{T}^{n}~,$ 
that is $\ell = n~, ~k = 1$ and the induction process 
yields the elements in $\cR_M (n)^n$ in this case. 
\end{thm}


\medbreak
So far, we have not made any assumptions which guarantee non--trivial
examples, but there is no doubt that there are many such situations, 
e.g. when $\rho: \pi_1 (B) \to \Aut (T) \cong \GL(d, \bZ)$ is 
surjective, $\rho: \pi_1 (B_g) \to \Aut (T^2) \cong \GL(2, \bZ)~,$ 
where $B_g$ is an oriented surface of genus $g>1~,$ 
or $\pi_1 (B)$ finite, etc. ...

\begin{rem}\label{degremark}
Degeneracy properties for the variety $\cS (n)$ 
and the representation variety $\cR_M (n)~:~$ 

\noindent
\begin{itemize}
\item[$(1)$]
The parametrization \eqref{generic3} of $\cS (n)$ is not 
closed in $\pi_1 (B) \bsl \what{T}_n~.$ At the limit points 
in \eqref{generic3} corresponding to $\cS (n)^{\ell}~, ~\ell < n~,$ 
the action of $\pi_1 (B)$ is still transitive, even though the orbit 
degenerates and the corresponding representations of $\pi_1 (M)$ are 
still irreducible. 
We denote by $\wbar{\cS (n)}$ the set of all 
$(m-d)$--dimensional submanifolds $\Sigma \hra \what M$ 
which are a finite union of leaves of the 
foliation $\hcF~,$ whose transversal holonomy is of order 
$\leq n$ (modulo transversal automorphisms of $\hcF$). 
In a precise sense, $\wbar{\cS (n)}$ corresponds to the closure 
of $\pi_1 (B) \bsl \what{T}^n$ in the orbit space 
$\pi_1 (B) \bsl  \what{T}_n$ of leaves with finite 
holonomy of order $\leq n~.$ 
The generic elements $\cS(n)^n$ correspond to the open, dense 
subset $\pi_1 (B) \bsl \what{T}^n \subset 
\wbar{\pi_1 (B) \bsl \what{T}^n} \subset 
\pi_1 (B) \bsl \what{T}_n~.$ 

\item[$(2)$]
At a limit point of $\pi_1 (B) \bsl \what{T}^{n}$ in 
$\pi_1 (B) \bsl \what{T}_{n-1} = 
\pi_1 (B) \bsl (\what{T}_{n} \bsl \what{T}^{n})~,$ 
the action of the holonomy group will generally fail to be 
transitive, the orbit structure degenerates, the index drops 
and we end up with a finite number of orbits, say $k_i$ times 
an orbit of order $\ell_i~,$ such that $\sum_i ~k_i \ell_i = n~.$ 
Geometrically, this means that the $n$--fold covering (leaf) 
$\hpis : \Sigma \to B$ collapses under this limiting process to 
$\ell_i$--fold coverings (leaves) $\hat\pi_{\Sigma_i}:\Sigma_i \to B$ 
of multiplicity $k_i~,$ satisfying the above relation. 
The `degeneracy' condition for 
$\Sigma \in \cS (n) \subset \wbar{\cS (n)}$ 
at the limit is then of the form 
\begin{equation}\label{generic5}
\begin{aligned}
(\underline{\bC}, \trivcon, \Sigma) &\bra 
\bigoplus_i ~k_i ~(\underline{\bC}, \trivcon, \Sigma_i)~, \\
k_i ~(\underline{\bC}, \trivcon, \Sigma_i) = 
(\underline{\bC}^{k_i}, \trivcon, \Sigma_i) &\in \Lochmn{\ell_i}~, ~
\sum_i ~k_i \ell_i = n~.
\end{aligned}
\end{equation}

\item[$(3)$] 
The limit degeneracy of the spectral covering $\Sigma$ corresponds 
of course to the degeneracy of the representations 
$\{\xi_1, \ldots, \xi_n\}$ of $\Lambda~.$ 
In fact, with multiple holonomy representations in the limit, the 
action of the holonomy group will generally fail to be transitive, 
the orbit structure will degenerate, the index will fall and the 
representation will decompose into sums of 
$k_i$ times an irreducible representations of rank 
$\ell_i~,$ such that $\sum_i ~k_i \ell_i = n~.$ 

\item[$(4)$] 
The space $\cR_M (n)$ of irreducible $\U (n)$--representations 
may be completed as well by adding (sums of) irreducible 
representations of lower rank, as described above. 
This is similar to the completion of stable bundles to include 
semistable bundles. 
Then $\Psi (n)$ extends by continuity to a surjective mapping 
\begin{equation}\label{generic6}
\wbar{\Psi (n)} : \wbar{\cR_M (n)} \lra \wbar{\cS (n)}~.
\end{equation}
Theorem \ref{mainresult4} $(1)$ remains valid on the generic 
(open dense) subset $\cR_M (n)^n~,$ but on the boundary of 
$\wbar{\cS (n)}$ the structure of the fibers of $\wbar{\Psi (n)}$ 
is more complicated, in accordance with Theorem \ref{mainresult4} 
$(2)$ and the previous remarks. 
\end{itemize}
\end{rem}

\begin{rem}\label{repremark}
Structure of the representation variety $\cR_M (n)~:~$ 
Theorem \ref{mainresult4} $(1)$ means that the representation 
variety of a torus bundle $\pi:M \to B$ resembles generically an 
{\it integrable system}, that is a fibration by abelian groups. 
It would be very interesting to determine the conditions under 
which $\cR_M (n)$ is a symplectic manifold, with the fibers of 
$\Psi (n) : \cR_M (n) \lra \cS(n)$ being Lagrangian over 
$\cS (n)^n~.$
\end{rem}



\subsection{Instantons on $T^1$--fibered $4$--manifolds} 
\label{instanton}

Here we consider the case $m=4~, ~d=1~,$ and take 
$g_M = g_{T(\pi)} \oplus \pi^* g_B$ to be a 
bundle--like Riemannian metric on $M$ with respect to the fiber 
space 
\eqref{torus.bundle1} and the exact sequence \eqref{exact1}. 
Assuming $M$ to be oriented, 
$g_M$ induces a splitting of the bundle of 2--forms on $M$ into 
selfdual (SD) and anti--selfdual (ASD) $2$--forms under 
the Hodge operator $\ast~:$ 
\begin{equation}\label{asd1}
\Omega^2_M = \Omega_M^+ \oplus \Omega_M^-~.
\end{equation}
From \eqref{bigrad1} we also have the decomposition 
\begin{equation}\label{bigrad2}
\Omega_M^2 \cong \Omega_M^{2,0} \oplus \Omega_M^{1,1} = 
\pi^* \Omega_B^2 \oplus \pi^* \Omega_B^1 \otimes \Omega_{M/B}^1~.
\end{equation}
Since $g_M$ is bundle--like, the Hodge operator exchanges the 
summands in \eqref{bigrad2} and therefore a $2$--form 
$\omega = (\omega^{2,0}, \omega^{1,1})$ satisfies 
$* \omega = \pm ~\omega$ if and only if 
\begin{equation}\label{asd1a}
* \omega = (* \omega^{1,1}, * \omega^{2,0}) = 
\pm~(\omega^{2,0}, \omega^{1,1}) = \pm~\omega~, 
\end{equation}
that is $* \omega^{1,1} = \pm ~\omega^{2,0}$ or equivalently 
$* \omega^{2,0} = \pm ~\omega^{1,1}$~, 
so that the projections $\Omega_M^2 \to \Omega_M^{2,0}$ and 
$\Omega_M^2 \to \Omega_M^{1,1}$ induce isomorphisms 
$$
\Omega_M^{\pm} \osetl{\cong} \Omega_M^{2,0} \qquad , \qquad 
\Omega_M^{\pm} \osetl{\cong} \Omega_M^{1,1}~. 
$$
Given a Hermitian vector bundle $E$ with unitary connection
$\nabla_E$ over $M$, recall that $\nabla_E$ is said to be SD
(resp. ASD) if its curvature $\nabla_E^2$ is SD (resp. ASD) 
as a $\End_s(E)$--valued 2--form, that is 
\begin{equation}\label{asd2}
\ast ~\nabla_E^2 = \pm ~\nabla_E^2~, 
\end{equation}
which by \eqref{asd1a} is equivalent to 
\begin{equation}\label{asd3}
(\nabla_E^2)^{2,0} = \pm~ \ast (\nabla_E^2)^{1,1}~. 
\end{equation}

\medbreak
Now $Z_{\Sigma}$ is $4$--dimensional and the metric
$\ps^*~g_M$ prescribes a Riemannian ramified covering
$(Z_{\Sigma},~ \ps^*~ g_M) \to (M, ~g_M)$~. 
The pull--back $\ps^*$ defines decompositions like 
\eqref{asd1} and \eqref{bigrad2}  on $Z_\Sigma$ relative 
to the pull--back fiber bundle $\hps : Z_\Sigma \to \what{M}$~. 
Using \eqref{curv1}, we see that the $\ASD$ equation \eqref{asd3}
on $Z_\Sigma$ can be written as  
\begin{equation}\label{asd4} \tilde{j}^* (\wti{\nabla}_E^2)^{2,0} = \ps^* (\nabla_E^2)^{2,0} = 
\pm~ \ps^* (\ast (\nabla_E^2)^{1,1}) = \pm~ \ast_\Sigma 
\ps^* (\nabla_E^2)^{1,1}~. 
\end{equation}
Suppose now that the adapted unitary connection $\nabla_E$ 
is in addition 
Poincar\'e basic. Then we have from \eqref{asd4} and  
\eqref{pbasic2} 
\begin{equation}\label{asd5}
\tilde{j}^* (\wti{\nabla}_E^2)^{2,0} = 
\ps^* (\nabla_E^2)^{2,0} = 
\pm~ \ast_\Sigma (\bF \mid {Z_\Sigma})~. 
\end{equation}

\medbreak
In the following Lemma 
we use the functor $\fm : \vectn \lra \specn$ to transform an 
instanton $(E, \nabla_E)$ on $M$ to the corresponding spectral 
data $(\cL,\nabla_{\cL},\Sigma)$~. 

\begin{lem} \label{asd6}
Suppose that the Poincar\'e basic unitary connection $\nabla_E$ 
is $\ASD$~. Then we have 
\begin{equation}\label{asd7}
\ps^* (\nabla_E^2)^{2,0} = 
\pm~ \ast_\Sigma (\bF \mid {Z_\Sigma})~. 
\end{equation}
Further, the scalar $(2,0)$--form 
$\hat\omega = \ast_\Sigma ~(\bF \mid {Z_\Sigma})$ 
is harmonic and in particular $\hps$--basic, that is 
$\hat\omega = \hps^* \omega~.$ 
The curvature $\nabla_{\cL}^2$ of the transformed connection 
$\nabla_{\cL}$ is then given by 
\begin{equation} \label{transf.asd1}
\nabla_{\cL}^2 = \pm~ R^1 \hpsstar \left( \hat\omega \right) = 
\pm ~\omega~. 
\end{equation}
\end{lem}

\begin{proof} 
We need to show that $\hat\omega$ is harmonic. 
Since $\tilde{j}^* (\wti{\nabla}_E^2)^{0,2} = 
\tilde{j}^* (\wti{\nabla}_E^2)^{1,1} = 0$ by assumption, 
we have $\tilde{j}^* \wti{\nabla}_E^2 = 
\tilde{j}^* (\wti{\nabla}_E^2)^{2,0} = \ps^* (\nabla_E^2)^{2,0}~.$ 
Computing traces and using \eqref{asd5}, we obtain 
\begin{equation} \label{chern1}
\tilde{j}^* \Tr \wti{\nabla}_E^2 = \pm~ n ~\ast_\Sigma 
(\bF \mid {Z_\Sigma})~. 
\end{equation}
Since the first Chern polynomial $\Tr \wti{\nabla}_E^2$ is closed, 
we see from \eqref{chern1} that 
\begin{equation} \label{chern2}
d \hat\omega = d \ast_\Sigma (\bF \mid {Z_\Sigma}) = 0~. 
\end{equation}
As $d \bF = d \nabla_{\rpb}^2 = 0$ from \eqref{p-curv}, 
it follows that $\bF \mid {Z_\Sigma}$ must be a {\em harmonic} 
$2$--form with respect to the bundle--like metric 
$g_\Sigma = \ps^*~ g_M$ on $Z_\Sigma~.$ 
Since $\hat\omega = \ast_\Sigma ~(\bF \mid {Z_\Sigma})$ is of type 
$(2,0)~,$ we have also $i_X \hat\omega = 0~,$ for any vector 
field $X$ in $T (\hps)$ and therefore 
$L_X \hat\omega = i_X d \hat\omega = 0~;$ in other words, 
$\hat\omega$ is in addition a $\hps$--basic form and so 
$\hat\omega = \hps^* \omega~,$ for a unique closed $2$--form $\omega$ 
on $\Sigma~.$ Equation \eqref{transf.asd1} follows then from 
\eqref{transf.curv2}. In fact, we have 
$\nabla_{\cL}^2 = \pm~ R^1 \hpsstar (\hat\omega) = 
\pm~ R^1 \hpsstar (\hps^* \omega) = \pm ~\omega~.$ 
\end{proof}

\medbreak
Note that the curvature term $\ps^* (\nabla_E^2)^{2,0}$ must be 
independent of the choice of the connection $\nabla_E~,$ since 
$\bF \mid Z_\Sigma$ depends only on the foliated structure 
$(E, \npart_E)~.$ Moreover, the curvature term 
$\ps^* (\nabla_E^2)^{2,0}$ must also be scalar--valued 
(i.e. assume values in the center of $\ps^* \End_s (E)$~), 
for $\bF$ is scalar--valued.

\begin{thm}\label{asdthm}
Assume that there exists a bundle--like metric $g_M$ with respect to which
the form  $\hat{\omega}=\ast_\Sigma ~(\bF \mid {Z_\Sigma})$ is harmonic.
The functors $\fm$ and $\fmh$ induce an equivalence between the following
objects$~:$
\begin{itemize}
\item[$(1)$]
Foliated Hermitian vector bundles $(E, \nabla_E) \in \vectn$ 
with Poincar\'e basic unitary connections $\nabla_E~,$ 
satisfying the $\ASD$--equation \eqref{asd3}. 

\item[$(2)$]
Relative skyscrapers $(\cS, \nabla_{\cS}, \Sigma)$ in $\specn~,$ 
such that the curvature $\nabla_{\cS}^2$ of the connection 
$\nabla_{\cS}$ satisfies 
\begin{equation}\label{transf.asd3}
\nabla_{\cS}^2 = \pm ~\omega~. 
\end{equation}
\end{itemize}
\end{thm}

The harmonicity condition \eqref{chern2} for the curvature 
$\bF \mid {Z_\Sigma}$ of the connection $\nabla_{\rpb}$ on 
$\rpb_\Sigma$ depends only on the foliated structure $(E, \npart_E)$ 
and the bundle--like metric $g_M$ and is therefore an a priori 
obstruction for the existence of Poincar\'e basic instantons, that is 
solutions of equation \eqref{asd3}, respectively \eqref{asd7}. 

\begin{proof}
This follows from combining Lemma \ref{asd6} with 
Theorem \ref{mainresult2}. 
\end{proof}

\medbreak
Finally, we analyze the properties of the parameter spaces for 
the various structures for a {\em fixed} foliated vector bundle 
$(E, \npart_E)~.$ For two adapted connections 
$\nabla_E~, \nabla_E'~,$ we have $\nabla_E' = \nabla_E + \vphi~,$ 
where $\vphi \in \Cinfty{M}{\End_s (E) \ot \Omega_M^{1,0}}~,$ 
that is the adapted connections form an affine space modeled on 
the linear space $\Cinfty{M}{\End_s (E) \ot \Omega_M^{1,0}}~.$ 
For $\nabla_E' = \nabla_E + \vphi~,$ we have also 
\begin{equation}\label{param1}
(\nabla^{' 2}_E)^{1,1} = (\nabla_E^2)^{1,1} + \npart_E (\vphi)~. 
\end{equation}
Therefore if $\nabla_E$ is {\em basic}, then the curvature term 
$(\nabla_E^2)^{1,1}$ vanishes and $\nabla_E' = \nabla_E + \vphi$ 
is basic if and only if $\npart_E (\vphi) = 0~.$ 
Thus the space of basic connections is either empty or else an 
affine space modeled on the linear space of $\npart$--parallel 
sections in $\End_s (E) \ot \Omega_M^{1,0}~.$ 
 
Now if $\nabla_E$ is {\em Poincar\'e basic}, then the curvature term 
$\tilde{j}^* (\wti{\nabla}_E^2)^{1,1} = 0$ and 
$\ps^* (\nabla_E^2)^{1,1}$ is fixed by \eqref{pbasic2}. 
Then $\nabla_E' = \nabla_E + \vphi$ is Poincar\'e basic, if and 
only if $\ps^* \npart_E (\vphi) = 0$ on $Z_\Sigma$ which is 
equivalent to $\npart_E (\vphi) = 0~.$ 
Thus the space of Poincar\'e basic connections on $E$ is an affine 
space modeled also on the linear space of $\npart$--parallel sections 
in $\End_s (E) \otimes \Omega_M^{1,0}~.$ 

For the instantons in this Section \ref{instanton}, 
the curvature term $\tilde{j}^* (\wti{\nabla}_E^2)^{2,0}$ 
is also fixed by \eqref{asd5} and $\nabla_E' = \nabla_E + \vphi$ 
satisfies the instanton equation \eqref{asd5}, if and only 
if in addition to the previous condition, the parameter 
$\vphi$ satisfies the quadratic PDE  
\begin{equation}\label{instparam1}
\nabla_E (\vphi) + \frac{1}{2} ~[\vphi, \vphi] = 0~. 
\end{equation}
Here the Lie bracket is taken in the adjoint bundle $\End_s (E)~.$ 
Note that the expressions in \eqref{instparam1}
are of type $(2,0)~,$ since we already have $\npart_E (\vphi) = 0$ 
from $(2)~.$

\subsection{Monopoles on $T^1$--fibered $3$--manifolds} 
\label{monopole1}

We keep the assumptions and notation of the previous Section 
\ref{instanton}, but now we take $m=3~, ~d=1~.$ 
The Hodge operator given by the bundle--like metric $g_M$ 
transforms now 
$$
\ast ~:~ \Omega_M^{2,0} \osetl{\cong} \Omega_M^{0,1}  \qquad , \qquad 
\ast ~:~ \Omega_M^{1,1} \osetl{\cong} \Omega_M^{1,0}~. 
$$
Given a Hermitian vector bundle $E$ with unitary connection
$\nabla_E$ over $M$, we consider the corresponding connection 
form $A$ on the unitary frame bundle $F_U (E)$ with curvature 
form $F_A~.$ 

\medbreak
The $\ASD$ equation is now replaced by the {\em monopole} 
equation relative to a Higgs field $\phi$ in $\End_s (E)~.$ 
\begin{equation}\label{mp1}
\ast ~F_A = D_A \phi = d \phi + [A, \phi]~, 
\end{equation}
or in terms of the corresponding unitary connection $\nabla_E$ 
\begin{equation}\label{mp2}
\ast ~\nabla_E^2 = \nabla_E (\phi)~. 
\end{equation}
The type decomposition of \eqref{mp1}~, respectively \eqref{mp2} 
is then given by
\begin{equation}\label{mp3}
\begin{aligned}
\ast ~F_A^{2,0} = D_A^{0,1} \phi &= d^{0,1} \phi + [A^{0,1}, \phi]~, \\
\ast ~F_A^{1,1} = D_A^{1,0} \phi &= d^{1,0} \phi + [A^{1,0}, \phi]~, 
\end{aligned}
\end{equation}
or in terms of the corresponding unitary connection $\nabla_E$ 
\begin{equation}\label{mp4}
\begin{aligned}
\ast ~(\nabla_E^2)^{2,0} &= \npart_E (\phi)~,   \\
\ast ~(\nabla_E^2)^{1,1} &= \nabla_E^{1,0} (\phi)~. 
\end{aligned}
\end{equation}

\medbreak
Now $Z_{\Sigma}$ is $3$--dimensional and the metric
$\ps^*~g_M$ prescribes a Riemannian ramified covering
$(Z_{\Sigma},~ \ps^*~ g_M) \to (M, ~g_M)$~. 
From \eqref{mp4} it follows that the monopole equations 
\eqref{mp4} for the pair $(A, \phi)~,$ 
respectively $(\nabla_E, \phi)$ on $M$ are now expressed on 
$Z_\Sigma$ by 
\begin{equation}\label{mp5}
\begin{aligned}
\ast_\Sigma ~\ps^* (\nabla_E^2)^{2,0} &= \ps^* \npart_E (\phi)~,   \\
\ast_\Sigma ~\ps^* (\nabla_E^2)^{1,1} &= \nabla^t_E (\ps^* \phi) = 
\ps^* \nabla_E^{1,0} (\phi)~. 
\end{aligned}
\end{equation}
In order to proceed with the reduction to $\Sigma~,$ we need to 
assume that the Higgs field $\phi$ is parallel along the fibers, 
that is $\npart_E (\phi) = 0~.$ By \eqref{mp4}, this is equivalent 
to $(\nabla_E^2)^{2,0} = 0~.$ 
Therefore we restrict attention to special solutions of the 
monopole equations \eqref{mp4}, namely 
\begin{equation}\label{mp6}
\begin{aligned}
(\nabla_E^2)^{2,0} = 0 \qquad &, \qquad \npart_E (\phi) = 0~,   \\
\ast ~(\nabla_E^2)^{1,1} &= \nabla_E^{1,0} (\phi)~. 
\end{aligned}
\end{equation}
Suppose again that the adapted unitary connection $\nabla_E$ is in 
addition Poincar\'e basic. Then we have 
$\tilde{j}^* (\wti{\nabla}_E^2)^{1,1} = 
\ps^* (\nabla_E^2)^{1,1} - \bF \mid Z_\Sigma = 0$ 
and therefore the connection $\tilde{j}^* \wti{\nabla}_E$ 
must be flat by Lemma \ref{anticomlemma}. 
On $Z_\Sigma~,$ the monopole equations \eqref{mp6} are now 
given by 
\begin{equation}\label{mp7}
\begin{aligned}
(\nabla_E^2)^{2,0} = 0 \qquad &, \qquad \npart_E (\phi) = 0~,   \\
\nabla_E^t (\ps^* \phi) = \ps^* \nabla_E^{1,0} (\phi) &= 
\ast_\Sigma ~(\bF \mid {Z_\Sigma})~. 
\end{aligned}
\end{equation}
Observe that the second equation is of type $(1,0)~.$ 
From \eqref{mp7} and the vanishing of the 
commutator $\Xi$ in \eqref{xisigma}, it follows that 
the scalar $(1,0)$--form $\hat\omega = \ast_\Sigma ~(\bF \mid {Z_\Sigma})$
is closed along the fibers, that is 
\begin{equation}\label{fiberclosed1}
d_{\npart} \hat\omega  = 
d_{\npart} (\ast_\Sigma ~(\bF \mid {Z_\Sigma})) = 0~. 
\end{equation}
From the first equation in \eqref{mp6} and \eqref{curv1}, we see 
that $(\nabla_E^t)^2 = 0$ and therefore \eqref{mp7} implies that 
\begin{equation}\label{transclosed1}
\nabla_E^t ~\hat\omega  = 
\nabla_E^t ~(\ast_\Sigma ~(\bF \mid {Z_\Sigma})) = 0~. 
\end{equation}
As $\hat\omega$ is scalar--valued, \eqref{fiberclosed1} and 
\eqref{transclosed1} are equivalent to 
\begin{equation}\label{closed1}
d ~\ast_\Sigma ~(\bF \mid {Z_\Sigma})) = 
d \hat\omega = \nabla_E ~\hat\omega = 0~.  
\end{equation}
It follows that the $(1,1)$--form $\bF \mid Z_\Sigma$ is 
{\em harmonic}, since $d \bF = 0$ from \eqref{p-curv}.  

\medbreak
As was the case for instantons in Section \ref{instanton}, 
the harmonicity of the curvature $\bF \mid Z_\Sigma$ 
of the connection $\nabla_{\rpb}$ on $\rpb_\Sigma$ 
is a necessary condition for the existence of solutions 
of \eqref{mp6}, depending only on the foliated structure 
$(E, \npart_E)$ and the bundle--like metric $g_M~.$ 

\medbreak
In the following Lemma we use again the functor 
$\fm : \vectn \lra \specn$ to transform a monopole 
$(E, \nabla_E, \phi)$ on $M$ to the corresponding 
spectral data $(\cL,\nabla_{\cL},\Sigma)$~. 

\begin{lem}\label{mp8}
Suppose that the Poincar\'e basic unitary connection $\nabla_E$ 
and the Higgs field $\phi$ satisfy the monopole equation \eqref{mp6}. 
Then the transformed connection $\nabla_{\cL}$ is flat, 
that is $\nabla_{\cL}^2 = 0~.$  
Further, the scalar $(1,0)$--form $\hat\omega = 
\ast_\Sigma ~(\bF \mid {Z_\Sigma})$ 
is harmonic and in particular $\hps$--basic, that is 
$\hat\omega = \hps^* \omega~.$ 
Setting $\bar\phi = R^1 \hpsstar (\phi_{\cK})~,$ the second 
equation in \eqref{mp7} transforms then into the equation  
\begin{equation}\label{transf.mp1}
\nabla_{\cL} \bar\phi = \nabla_{\cL} R^1 \hpsstar (\phi_{\cK}) = 
R^1 \hpsstar (\hat\omega) = \omega~. 
\end{equation}
\end{lem}

\begin{proof}
Recall that $\nabla_{\cL}$ is defined in \eqref{transf.conn} by 
$$
\nabla_{\cL} = \cR^1 \hpsstar (\nabla_E^{\ker})~:~
\cL \lra \cL \otimes \Omega^1_\Sigma~, 
$$
where $\cL = \cR^1 \hpsstar (\cK)~.$ 
The flatness of $\nabla_{\cL}$ follows from \eqref{transf.curv2}, 
since $(\nabla_E^2)^{2,0} = 0~.$ From $\npart_E (\phi) = 0$ 
it follows that $\ps^* \phi$ preserves 
the sheaf $\cK~.$ Thus $\phi$ induces an endomorphism 
$\phi_{\cK} = \ps^* \phi\mid \cK : \cK \to \cK$ 
and $\nabla_E^t$ determines the homomorphism 
$\nabla_E^{\ker} \phi_{\cK} : \cK \to \cK \otimes \hps^* \Omega_\Sigma^1~,$ 
by the usual formula 
$\nabla_E^{\ker} (\phi_{\cK}) = \nabla_E^{\ker} \circ \phi_{\cK} - 
(\phi_{\cK} \otimes \id) \circ \nabla_E^{\ker}~.$ 
Using \eqref{fiberclosed1}, we may now rewrite the second equation in 
\eqref{mp7} as $\nabla_E^{\ker} (\phi_{\cK}) = \hat\omega~.$ 
Thus by \eqref{transf.conn}, this transforms into the equation 
\begin{equation}\label{mp9}
\nabla_{\cL} R^1 \hpsstar (\phi_{\cK}) = 
R^1 \hpsstar (\hat\omega)~. 
\end{equation}
As already noted, the form $\hat\omega$ is harmonic by 
\eqref{closed1}. Since $\hat\omega$ is of type 
$(1,0)~,$ we have also $i_X \hat\omega = 0~,$ for any vector 
field $X$ in $T (\hps)$ and therefore from \eqref{fiberclosed1} 
$L_X \hat\omega = i_X d \hat\omega = i_X d_{\npart} \hat\omega = 0~;$ 
in other words, $\hat\omega$ is in addition a $\hps$--basic form and so 
$\hat\omega = \hps^* \omega~,$ for a unique closed $2$--form $\omega$ 
on $\Sigma~.$ Hence 
$\nabla_{\cL} R^1 \hpsstar (\phi_{\cK})  = \omega$ from \eqref{mp9}.  
The fact that $\omega~,$ hence $\hat\omega = \hps^* \omega$ are 
closed follows of course also from the flatness of $\nabla_{\cL}~,$ 
since we have $0 = \nabla_{\cL}^2 R^1 \hpsstar (\phi_{\cK}) = d \omega~.$ 
\end{proof}

\begin{thm}\label{monopolthm1}
The functors $\fm$ and $\fmh$ induce an equivalence between the
following objects$~:$
\begin{itemize}
\item[$(1)$]
Foliated Hermitian vector bundles $(E, \nabla_E) \in \vectn$ 
with Poincar\'e basic unitary connections $\nabla_E$ and 
$\npart$--parallel Higgs fields $\phi$ in $\End_s (E)~,$ 
satisfying the monopole equations \eqref{mp6}~. 

\item[$(2)$]
Relative skyscrapers $(\cS, \nabla_{\cS}, \Sigma) \in \specn~,$ 
such that $\nabla_{\cS}$ is flat, and Higgs fields 
$\phi_{\cS}$ in $\End (\cS)~,$ satisfying 
\begin{equation}\label{transf.mp2}
\nabla_{\cS} \phi_{\cS} = \omega~. 
\end{equation}
\end{itemize}
\end{thm}

\begin{proof}
This follows from Lemma \ref{mp8} combined with 
Theorem \ref{mainresult2}. 
\end{proof}

\medbreak
To conclude, let us analyze the properties of the parameter space 
of monopoles $(\nabla_E, \phi)$ for a {\em fixed} foliated 
vector bundle $(E, \npart_E)~.$ 
For two Poincar\'e basic connections $\nabla_E~, \nabla_E'~,$ 
we have again $\nabla_E' = \nabla_E + \vphi~,$ where 
$\vphi \in \Cinfty{M}{\End_s (E) \ot \Omega_M^{1,0}}~,$ 
satisfies $\npart_E (\vphi) = 0$ as in the remarks at 
the end of Section \ref{instanton}. 
For the monopoles in this Section \ref{monopole1}, 
the curvature term $\tilde{j}^* (\wti{\nabla}_E^2) = 0~.$ 
Hence $\nabla_E' = \nabla_E + \vphi$ satisfies the monopole 
equations \eqref{mp6}, if and only if $\vphi$ satisfies the 
quadratic PDE  
\begin{equation}\label{mpparam1}
\nabla_E (\vphi) + \frac{1}{2} ~[\vphi, \vphi] = 0~. 
\end{equation} 
Note that the expressions in \eqref{mpparam1} 
are of type $(2,0)~,$ since we already have 
$\npart_E (\vphi) = 0~.$  
Thus the parameter space for the monopoles is the same as 
for the instantons (compare \eqref{instparam1}). 
The monopole equations \eqref{mp6} are linear 
in the Higgs fields $\phi~.$ Therefore $\phi'$ also satisfies 
\eqref{mp6}, respectively \eqref{mp7}, if and only if 
$\phi' = \phi + \psi~,$ where $\psi \in \Cinfty{M}{\End_s (E)}$ 
satisfies 
\begin{equation}\label{mpparam2}
\nabla_E (\psi) = 0~. 
\end{equation}




\section*{Appendix : DeRham complexes along the fibers} \label{append1}

Throughout this paper we made use of the fiberwise DeRham complex 
relative to a fiber bundle $\pi : M \to B~.$ 
This complex is well known from foliation theory; 
in our context is extensively used in \cite{BMP1, BMP2}. 

\medbreak
For any foliated vector bundle $(E, \npart_E)$ on $M~,$ 
there is a fiberwise DeRham complex of sheaves 
\begin{equation}\label{leafres1}
0 \to \ker \npart_{E} \otimes \pi^* \Omega_B^u 
\osetl{\varepsilon} E \otimes \Omega_M^{u,0} \osetl{d_{\npart}} 
E \otimes \Omega_M^{u,1} \osetl{d_{\npart}} 
E \otimes \Omega_M^{u,2} \osetl{d_{\npart}} \ldots 
\osetl{d_{\npart}} E \otimes \Omega_M^{u,d}~, 
\end{equation}
which is a fine resolution of the sheaf of $\npart$--parallel 
sections in $E \otimes \Omega_M^{u,0}~.$ 
Therefore the sheaves $E \otimes \Omega_M^{u,*}$ 
are $p_*$--acyclic and also $\Gamma (M, )$--acyclic for the global 
section functor $\Gamma (M, )~,$ that is the derived 
direct images $\cR^j \pi_* (E \otimes \Omega_M^{u,*}) = 0~, ~j > 0~.$  
It follows that the higher direct images of 
$\ker(\npart_{E}) \otimes \pi^* \Omega_B^u$ can be computed 
from the fine resolution \eqref{leafres1} and the projection 
formula by 
\begin{equation}\label{deriv1}
\cR^j \pi_* (\ker \npart_{E} \otimes \pi^* \Omega_B^u) \cong 
\cR^j \pi_* (\ker \npart_{E}) \otimes \Omega_B^u \cong 
\cH^j (\pi_* (E \otimes \Omega^{0,*}, d_{\npart})) 
\otimes \Omega_B^u ~. 
\end{equation}
Likewise, the global fiberwise cohomology is given by 
\begin{equation}\label{deriv2}
H^{u,j}_{\pi} (M, E) \equiv 
H^j (M, \ker \npart_{E} \otimes \pi^* \Omega_B^u) \cong 
H^j (\Gamma(M, E \otimes \Omega^{u,*}), d_{\npart})~.  
\end{equation}
The two cohomologies are linked by the convergent Leray 
spectral sequence 
\begin{equation}\label{specseq1}
E_2^{i,j} = H^i (B, \cR^j \pi_* (\ker \npart_{E}) \otimes \Omega_B^u) 
\Ra H^{i+j} (M, \ker \npart_{E} \otimes \pi^* \Omega_B^u)~, 
\end{equation}
with edge homomorphisms
\begin{equation}\label{edge1}
E_2^{j,0} = H^j (B, \pi_* (\ker \npart_{E}) \otimes \Omega_B^u) \ra  
H^j (M, \ker \npart_{E} \otimes \pi^* \Omega_B^u) \ra 
E_2^{0,j} = R^j \pi_* (\ker \npart_{E} \otimes \pi^* \Omega_B^u)~,  
\end{equation}
where we set $R^j \pi_* ( \cdot ) = \Gamma (B, \cR^j \pi_* ( \cdot ))$
for the global sections in $\cR^j \pi_* ( \cdot )~.$
In our context of torus fiber bundles, we encounter vanishing conditions, 
leading to degeneracy conditions for the spectral sequence. 
If $\cR^j \pi_* (\ker \npart_{E}) = 0~, ~0 < j \leq d~,$ 
the non--zero terms are determined by edge isomorphisms 
\begin{equation}\label{degen1}
E_2^{j,0} = H^j (B, \pi_* (\ker \npart_{E}) \otimes \Omega_B^u) 
\osetl{\cong} 
H^{j} (M, \ker \npart_{E} \otimes \pi^* \Omega_B^u)~, ~j \geq 0~. 
\end{equation}
If $\cR^j \pi_* (\ker \npart_{E}) = 0~, ~0 \leq j < d~,$ 
the non--zero terms are determined by edge isomorphisms 
\begin{equation}\label{degen2}
H^{d+j} (M, \ker \npart_{E} \otimes \pi^* \Omega_B^u) \osetl{\cong} 
E_2^{j,d} = H^j (B, \cR^d \pi_* (\ker \npart_{E}) \otimes \Omega_B^u)~, ~
j \geq 0~.  
\end{equation}
In particular, we have for $j=0$~:
\begin{equation}\label{degen3}
H^{d} (M, \ker \npart_{E} \otimes \pi^* \Omega_B^u) \osetl{\cong} 
E_2^{0,d} = \Gamma (B, \cR^d \pi_* (\ker \npart_{E}) \otimes \Omega_B^u) = 
R^d \pi_* (\ker \npart_{E} \otimes \pi^* \Omega_B^u)~. 
\end{equation}

\medbreak
The previous discussion of basic connections in Section \ref{basic1} 
could have been formulated in terms of this fiberwise resolution 
(and its global cohomology) with coefficients in the foliated adjoint 
vector bundle $\End_s (E)~.$ Specifically, the mixed curvature term 
$(\nabla_E^2)^{1,1}$ of an adapted connection $\nabla_E$ for 
$(E, \npart_E)$ satisfies $d_{\npart} (\nabla_E^2)^{1,1} = 0$ 
and for $\nabla_E' = \nabla_E + \vphi~, ~\vphi \in 
\Gamma(M, \End_s (E) \otimes \Omega_M^{1,0})~,$ we have from 
$$
(\nabla^{' 2}_E)^{1,1} = (\nabla_E^2)^{1,1} + \npart_E (\vphi) = 
(\nabla_E^2)^{1,1} + d_{\npart} \vphi~. 
$$
Thus $(\nabla_E^2)^{1,1}$ defines a cohomology class 
\begin{equation}\label{atiyah1}
a (E, \npart_E) = [(\nabla^2_E)^{1,1}] \in 
H^{1,1}_{\pi} (M, \End_s (E)) = 
H^1 (\Gamma(M, \End_s (E) \otimes \Omega^{1,*}, d_{\npart}))~, 
\end{equation}
depending only on the foliated vector bundle $(E, \npart_E)~.$ 
This class is very similar to the Atiyah class in the 
theory of holomorphic vector bundles, where it obstructs the 
existence of a complex analytic connection. 
By construction, the class $a (E, \npart_E)$ is exactly 
the obstruction to the existence of a basic connection for 
$(E, \npart_E)~.$

In Sections \ref{fm1} and \ref{fm2}, 
the resolution \eqref{leafres1} is implicitly used with respect 
to the pull--back fiber bundle $\hat{p} : Z \to \what{M}~,$ 
the fiberwise derivative $\nabla_E^r$ and its restriction to 
$\hps : Z_\Sigma \to \Sigma~.$


\noindent{\bf Acknowledgment}
J. F. G. and F. W. K. gratefully acknowledge the hospitality
and support of the Erwin Schr\"odinger International Institute
for Mathematical Physics in Vienna (AT).  We also thank
Professors K. Burns and D. Gallo for various comments.

\end{document}